\documentclass{article}

\usepackage[a4paper, total={17cm, 24cm}]{geometry}
\usepackage[utf8]{inputenc}
\usepackage{amsthm}
\usepackage{amsmath}
\usepackage{amssymb,authblk}
\usepackage{amsrefs}
\usepackage[shortlabels]{enumitem}
\usepackage{tikz}
\usetikzlibrary{shapes}
\usetikzlibrary{positioning}
\usetikzlibrary{patterns}
\usetikzlibrary {decorations.pathmorphing}
\usepackage{xcolor}\usepackage{hyperref}
\usepackage{subcaption}
\numberwithin{equation}{section}

\newtheorem{thm}{Theorem}

\newtheorem{lem}[thm]{Lemma}

\newtheorem{defn}{Definition}

\newtheorem{conj}{Conjecture}
\newtheorem{claim}[thm]{Claim}

\title
{Minimally tough series--parallel graphs with \\ toughness at least $1/2$}

	\author[1,2]{Gyula Y. Katona}
	\author[1]{Humara Khan}
	\affil[1]{Department of Computer Science and Information Theory, Budapest University of Technology and Economics, Budapest, Hungary}

    \affil[2]{HUN-REN–ELTE Numerical Analysis and Large Networks Research Group, Budapest, Hungary}
	
\begin{document}
 
 	\maketitle
 
 \begin{abstract}
 	Let $t$ be a positive real number. A graph is called \emph{$t$-tough} if the removal of any vertex set $S$ that disconnects the graph leaves at most $|S|/t$ components. The toughness of a graph is the largest $t$ for which the graph is $t$-tough. A graph is minimally $t$-tough if the toughness of the graph is $t$, and the deletion of any edge from the graph decreases the toughness. Series--parallel graphs are graphs with two distinguished vertices called
 	terminals, formed recursively by two simple composition operations, series
 	and parallel joins. They can be used to model series and parallel electric circuits.
 	
 	We characterize minimally $t$-tough  series--parallel graphs for all $t\ge 1/2$. It is clear that there is no minimally $t$-tough  series--parallel graph if $t>1$. We show that, for $1\ge t >1/2$, most of the series--parallel graphs with toughness $t$ are minimally $t$-tough, but most of the series--parallel graphs with toughness $1/2$ are not minimally $1/2$-tough.

  \end{abstract}

	\section{Introduction}
	
	All graphs considered in this paper are undirected, but they may contain loops
	and parallel edges. Let $c(G)$ denote the number of components, and $\kappa(G)$
	the vertex-connectivity of $G$. For a connected graph $G$, a vertex set
	$S\subseteq V(G)$ is called a {\em cutset} if $c (G-S)>1$. The \emph{union} of
	two graphs $G_{1}=(V_{1},E_{1})$ and $G_{2}=(V_{2},E_{2})$ is the graph $G_{1}\cup
	G_{2}=(V_{1}\cup V_{2}, E_{1}\cup E_{2})$. In the present paper, we only use
	this notation when $E_{1}\cap E_{2}=\emptyset$ and $|V_{1}\cap V_{2}|\le 2$. For
	an edge $uv$ where $u,v\in V_{1}$ the union of $G_{1}$ and $uv$ is $G_{1}\cup
	uv=(V_{1},E_{1}\cup uv)$.
	
	The notion of toughness was introduced by Chv\'{a}tal in \cite{toughness_intro}.
	
	\begin{defn}
		Let $t$ be a real number. A graph $G$ is called {\em $t$-tough} if
		$|S| \ge t \cdot c (G-S)$ for any set $S \subseteq V(G)$ with $c (G-S)>1$. The
		{\em toughness} of $G$, denoted by $\tau(G)$, is the largest $t$ for which G
		is $t$-tough, taking $\tau(K_{n}) = \infty$ for all $n \ge 1$. We say that a
		cutset $S \subseteq V(G)$ is a {\em tough set} if $|S|=\tau(G)\cdot c(G-S)$.
	\end{defn}
	Note that a graph is disconnected if and only if its toughness is $0$.
	
	Toughness is a widely used measure in network reliability and appears in applications involving various types of networks. So it is a natural
	question to investigate the following concept of minimally $t$-tough graphs, introduced
	by Broersma in \cite{t3}.
	
	\begin{defn}
		Let $t$ be a real number. A graph $G$ is said to be
		{\em minimally $t$-tough} if $\tau(G) = t$ and $\tau(G - e) < t$ for every
		$e \in E(G)$.
	\end{defn}
	
	It follows directly from the definition that every $t$-tough noncomplete graph
	is $2t$-connected, implying $\kappa(G) \ge 2 \tau(G)$ for noncomplete graphs.
	Therefore, the minimum degree of any $t$-tough noncomplete graph is at least $\lceil
	2t \rceil$.
	
	The following conjecture is motivated by a theorem of Mader \cite{mader},
	which states that every minimally $k$-connected graph has a vertex of degree
	$k$.
	
	\begin{conj}
		[Kriesell \cite{kriesell}] \label{kriesell1} Every minimally $1$-tough graph
		has a vertex of degree $2$.
	\end{conj}
	
	This conjecture can be naturally generalized.
	
	\begin{conj}
		[Generalized Kriesell's Conjecture \cite{specgraph}] \label{kriesell} Every minimally
		\allowbreak $t$-tough graph has a vertex of degree $\lceil 2t \rceil$.
	\end{conj}
	
	Recently, Conjecture~\ref{kriesell} was disproved. In \cite{zheng-sun} Zheng
	and Sun constructed minimally $\left(1+\frac{1}{2k-1}\right)$-tough 4-regular
	graphs for every $k\ge 2$. Since $\left\lceil 2\left(1+\frac{1}{2k-1}\right)\right
	\rceil<3$, the Conjecture~\ref{kriesell} is not true for all values of $t$, but
	it still might be true for all other values of $t$, or for all irregular
	graphs \cite{zheng-sun}.
	
	On the other hand, Conjecture~\ref{kriesell} was verified for several graph classes
	and various values of $t$. In \cite{specgraph} Conjecture~\ref{kriesell} was
	proved for minimally $t$-tough split graphs for all $t$, minimally $t$-tough
	chordal graphs when $t\leq1$, and minimally $t$-tough claw-free graphs when $t
	\leq 1$. Also, in \cite{t4} Ma et al.~showed that for minimally $\frac{3}{2}$-tough,
	claw-free graphs the conjecture is true.
	
	This paper focuses on series--parallel graphs, which are important in graph theory and computer science due to their structured and hierarchical properties.
	This inherent organization simplifies their analysis compared to general
	graphs. Moreover, series--parallel graphs are particularly valuable in
	applications such as network design, circuit layout, and algorithm development,
	as their structure allows many computationally intensive algorithms to run more
	efficiently.
	
	There are several ways to define series--parallel graphs. The following definition is based on that given in \cite{kawano}.
	
	\begin{defn}
		[Series--Parallel Graph] A graph $G(s,t)$ is a \emph{series--parallel graph}
		(sp-graph for short) with terminals $s$ and $t$, if either $G$ consists of one
		edge connecting $s$ and $t$, or $G$ is derived from two or more series--parallel
		graphs by one of the following two operations.
		\begin{itemize}
			\item \textit{Series join:} Given $k$ series--parallel graphs $G_{1}(s_{1},
			t_{1}),\allowbreak G_{2}(s_{2}, t_{2}), \allowbreak \ldots,\allowbreak G_{k}
			(s_{k}, t_{k})$, form a new graph $G(s, t)$ by identifying the vertices
			$t_{i}$ and $s_{i+1}$ for all $1 \leq i \leq k-1$ and setting $s=s_{1}$
			and $t=t_{k}$.
			
			\item \textit{Parallel join:} Given $k$ series--parallel graphs
			$G_{1}(s_{1}, t_{1}), G_{2}(s_{2}, t_{2}), \dots, G_{k}(s_{k}, t_{k})$, construct
			a new graph $G(s, t)$ by identifying the vertices $s_{1}, s_{2}, \dots, s
			_{k}$ as a single vertex $s$, and similarly identifying the vertices
			$t_{1}, t_{2}, \dots, t_{k}$ as a single vertex $t$.
		\end{itemize}
	\end{defn}
	
\begin{figure}[h]
	\begin{center}
			\scalebox{.9}{\begin{tikzpicture}[scale=.75, bend angle=20, line width=.7mm,	every edge/.style = {draw,line width=3mm}]
	
\node[draw, circle, fill=black, inner sep=1.5pt, label=above:$s$] (s) at (0,0) {};
\node[draw, circle, fill=black, inner sep=1.5pt, label=above:$s_1$] (s1) at (2,0) {};
\node[draw, circle, fill=black, inner sep=1.5pt, label=above:] (v1) at (4,2) {};
\node[draw, circle, fill=black, inner sep=1.5pt, label=below:] (v4) at (3.5,0) {};
\node[draw, circle, fill=black, inner sep=1.5pt, label=below:] (v5) at (4.5,0) {};
\node[draw, circle, fill=black, inner sep=1.5pt, label=below:] (v6) at (4,-2) {};
\node[draw, circle, fill=black, inner sep=1.5pt, label=above:$s_2$] (s2) at (6,0) {};
\node[draw, circle, fill=black, inner sep=1.5pt, label=above:] (x1) at (8,2) {};
\node[draw, circle, fill=black, inner sep=1.5pt, label=above:] (v7) at (11,0.5) {};
\node[draw, circle, fill=black, inner sep=1.5pt, label=above:] (v8) at (10,1.5) {};
\node[draw, circle, fill=black, inner sep=1.5pt, label=above:] (v9) at (12,1) {};
\node[draw, circle, fill=black, inner sep=1.5pt, label=below:] (x2) at (7.5,-1.5) {};
\node[draw, circle, fill=black, inner sep=1.5pt, label=below:] (x3) at (9,-3) {};
\node[draw, circle, fill=black, inner sep=1.5pt, label=below:] (v10) at (11,-1) {};
\node[draw, circle, fill=black, inner sep=1.5pt, label=below:] (v20) at (12,-2) {};
\node[draw, circle, fill=black, inner sep=1.5pt, label=above:$s_3$] (s3) at (14,0) {};
\node[draw, circle, fill=black, inner sep=1.5pt, label=above:] (v11) at (16,1.5) {};
\node[draw, circle, fill=black, inner sep=1.5pt, label=above:] (v12) at (16.5,0.5) {};
\node[draw, circle, fill=black, inner sep=1.5pt, label=above:] (v16) at (18,2) {};
\node[draw, circle, fill=black, inner sep=1.5pt, label=above:] (x7) at (21,2.5) {};
\node[draw, circle, fill=black, inner sep=1.5pt, label=above:] (v13) at (23,2) {};
\node[draw, circle, fill=black, inner sep=1.5pt, label=above:] (x8) at (19.5,1) {};
\node[draw, circle, fill=black, inner sep=1.5pt, label=above:] (v14) at (21,.3) {};
\node[draw, circle, fill=black, inner sep=1.5pt, label=below:] (x9) at (19,-2) {};
\node[draw, circle, fill=black, inner sep=1.5pt, label=above:] (v15) at (22,-2) {};
\node[draw, circle, fill=black, inner sep=1.5pt, label=right:$t$] (t) at (24,0) {};
\draw[] (s) to node[midway, above]{$e_{1}$} (s1);
\draw[] (s1) to node[midway, above]{$e_{2}$} (v1);
\draw[] (v1) to node[midway, above]{$e_{3}$} (s2);
\draw[] (s1) to node[midway, above]{$e_{4}$} (v4);
\draw[] (v4) to node[midway, above]{$e_{5}$} (v5);
\draw[] (v5) to node[midway, above]{$e_{6}$} (s2);
\draw[] (s1) to node[midway, below]{$e_{7}$} (v6);
\draw[] (v6) to node[midway, below]{$e_{8}$} (s2);
\draw[] (s2) to node[midway, above]{$e_{9}$} (x1);
\draw (x1) to node[midway, below]{$e_{10}$} (v7);
\draw (v7) to node[midway, above]{$e_{11}$} (s3);
\draw[bend left=50] (x1) to node[midway, above]{$e_{15}$} (s3);
\draw[] (x1) to node[midway, above]{$e_{12}$} (v8);
\draw[] (v8) to node[midway, above]{$e_{13}$} (v9);
\draw[] (v9) to node[midway, above]{$e_{14}$} (s3);
\draw[] (s2) to node[midway, left]{$e_{16}$} (x2);
\draw[] (x3) to node[midway, above]{$e_{18}$} (v10);
\draw[] (x2) to node[midway, left]{$e_{17}$} (x3);
\draw[] (v10) to node[midway, above]{$e_{19}$} (s3);
\draw[] (x3) to node[midway, below]{$e_{20}$} (v20);
\draw[] (v20) to node[midway, right]{$e_{21}$} (s3);
\draw[] (s3) to node[midway, above]{$e_{22}$} (v11);
\draw[] (v11) to node[midway, above]{$e_{23}$} (v16);
\draw[] (s3) to node[midway, below]{$e_{24}$} (v12);
\draw[] (v12) to node[midway, right]{$e_{25}$} (v16);
\draw[] (v16) to node[midway, above]{$e_{26}$} (x7);
\draw[] (x7) to node[midway, above]{$e_{27}$} (v13);
\draw[] (x7) to node[midway, left]{$e_{29}$} (t);
\draw[] (v13) to node[midway, right]{$e_{28}$} (t);
\draw[bend left=30] (v16) to node[midway, right]{$e_{30}$} (v14);
\draw[] (v16) to node[midway, below]{$e_{31}$} (x8);
\draw[] (x8) to node[midway, below]{$e_{32}$} (v14);
\draw[] (v14) to node[midway, above]{$e_{33}$} (t);
\draw[] (s3) to node[midway, below]{$e_{35}$} (x9);
\draw[] (x9) to node[midway, below]{$e_{37}$} (v15);
\draw[] (x9) to node[midway, above]{$e_{36}$} (t);
\draw[] (v15) to node[midway, right]{$e_{38}$} (t);
\draw[bend right=10] (s3) to node[midway, above]{$e_{34}$} (t);
\end{tikzpicture}}
		\caption{An example of a series--parallel graph with terminals $s$ and $t$.}\label{sp-graph}
	\end{center}
	\end{figure}
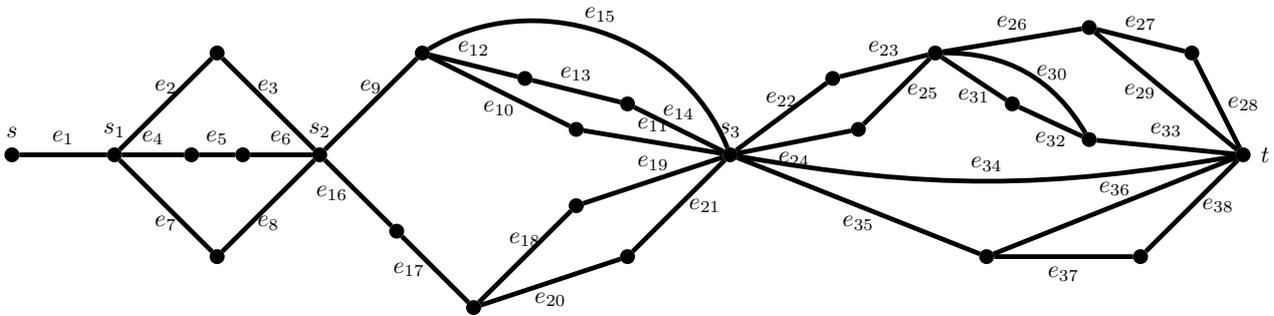

	{In Fig.~\ref{sp-graph} the subgraph formed by edges $e_1,\ldots, e_8$ is joined in series  with the subgraph formed by edges $e_9,\ldots, e_{21}$.  The subgraph formed by edges $e_{9},\ldots, e_{15}$ is joined in parallel with the subgraph formed by edges $e_{16},\ldots, e_{21}$. }
	
	In either operation, instead of joining $k$ graphs in one step, we could make
	$k-1$ joins each time joining precisely $2$ graphs.
		Note that the ordering
	$G_{1}(s_{1}, t_{1}), \allowbreak G_{2}(s_{2}, t_{2}), \allowbreak\ldots, \allowbreak
	G_{k}(s_{k}, t_{k})$
	matters for the series join, while it does not matter for the parallel join.
	
	The recursive definition of the series--parallel graph above naturally gives a
	rooted, leveled, vertex--ordered tree, called a \emph{series--parallel tree} (sp-tree
	for short), for each series--parallel graph $G(s, t)$. 
	To avoid confusion, we refer to the vertices of this tree as \emph{nodes},
	while the word \emph{vertex} will always refer to a vertex of the series--parallel graph.
	
	Each leaf in an sp-tree corresponds to an edge of $G(s, t)$, and each non-leaf
	node in an sp-tree corresponds to either a series or a parallel join. We call such
	nodes \emph{normal, series}, or \emph{parallel}, respectively. {The sp-tree of the graph of Fig.~\ref{sp-graph} is shown in  Fig.~\ref{sp-tree}. Each leaf in the sp-tree labeled $e_i$ represents the edge with the same label in Fig.~\ref{sp-graph}. The types of the internal nodes are also shown in  Fig.~\ref{sp-tree}.}
	
	The parent of a normal node can be either a series or a parallel node. The
	parent of a series or a parallel node could be either a series or a parallel node. However,
	it is easy to see that the operations can always be chosen so that the parent
	of a series node is always a parallel node, and the parent of a parallel node is
	always a series node. We can also assume that the children of series and
	parallel nodes are ordered from left to right in the same order in which the
	corresponding operations were performed. We denote this unique series--parallel
	tree of $G$ by $T_{G}$. {Fig.~\ref{sp-tree} shows this sp-tree.}
	
	\begin{figure}
		\begin{tikzpicture}[
every edge/.style = {draw,very thick},
vertex/.style args = {#1 #2}{circle, 
	draw, 
	thick,
	minimum size=2pt,inner sep=-2pt,
	fill=black,
	label=#1:#2},
leaf/.style={rectangle, draw=black, rounded corners=1mm, 
	text centered, anchor=north, text=black,text width=1em,inner sep=2.5pt},
parallel/.style={circle, draw=black, fill=gray!30,
	text centered, anchor=north, text=black,text width=1em,inner sep=1pt},
non/.style={circle, 
	text centered, anchor=north, text=black,text width=1em,inner sep=1pt},
series/.style={circle, draw=black, fill=black!70, 
	text centered, anchor=north, text=white, text width=1em,inner sep=1pt},
level distance=0.3cm, growth parent anchor=south]
\begin{scope}
\node[series] {$s$}
  child [sibling distance=25mm] {node [leaf] {$e_{1}$}}
  child [sibling distance=50mm] {node [parallel] {$p$}
    child [sibling distance=15mm] {node [series] {$s$}
      child [sibling distance=6mm] {node [leaf] {$e_{2}$}}
      child [sibling distance=6mm] {node [leaf] {$e_{3}$}}
    }
    child [sibling distance=15mm] {node [series] {$s$}
      child [sibling distance=6mm] {node [leaf] {$e_{4}$}}
      child [sibling distance=6mm] {node [leaf] {$e_{5}$}}
      child [sibling distance=6mm] {node [leaf] {$e_{6}$}}
    }
    child [sibling distance=15mm] {node [series] {$s$}
      child [sibling distance=6mm] {node [leaf] {$e_{7}$}}
      child [sibling distance=6mm] {node [leaf] {$e_{8}$}}
    }
  }
  child [sibling distance=40mm] {node [parallel] {$p$}
    child [sibling distance=20mm] {node [series] {$s$}
      child [sibling distance=6mm] {node [leaf] {$e_{9}$}}
      child [sibling distance=10mm] {node [parallel] {$p$}
        child [sibling distance=15mm] {node [series] {$s$}
          child [sibling distance=6mm] {node [leaf] {$e_{10}$}}
          child [sibling distance=6mm] {node [leaf] {$e_{11}$}}
        }
        child [sibling distance=20mm] {node [series] {$s$}
          child [sibling distance=6mm] {node [leaf] {$e_{12}$}}
          child [sibling distance=6mm] {node [leaf] {$e_{13}$}}
          child [sibling distance=6mm] {node [leaf] {$e_{14}$}}
        }
        child [sibling distance=10mm] {node [leaf] {$e_{15}$}}
      }
    }
    child [sibling distance=25mm] {node [series] {$s$}
    	child [sibling distance=6mm] {node [leaf] {$e_{16}$}}
    	child [sibling distance=6mm] {node [leaf] {$e_{17}$}}
      child [sibling distance=10mm] {node [parallel] {$p$}
        child [sibling distance=12mm] {node [series] {$s$}
        	child [sibling distance=6mm] {node [leaf] {$e_{18}$}}
        	child [sibling distance=6mm] {node [leaf] {$e_{19}$}}
        }
        child [sibling distance=12mm] {node [series] {$s$}
          child [sibling distance=6mm] {node [leaf] {$e_{20}$}}
          child [sibling distance=6mm] {node [leaf] {$e_{21}$}}
        }
      }
    }
  }
  child [sibling distance=65mm] {node [parallel] {$p$}
    child [sibling distance=15mm] {node [series] {$s$}
      child [sibling distance=30mm] {node [parallel] {$p$}
        child [sibling distance=12mm] {node [series] {$s$}
          child [sibling distance=6mm] {node [leaf] {$e_{22}$}}
          child [sibling distance=6mm] {node [leaf] {$e_{23}$}}
        }
        child [sibling distance=12mm] {node [series] {$s$}
          child [sibling distance=6mm] {node [leaf] {$e_{24}$}}
          child [sibling distance=6mm] {node [leaf] {$e_{25}$}}
        }
      }
      child [sibling distance=30mm] {node [parallel] {$p$}
        child [sibling distance=22mm] {node [series] {$s$}
          child [sibling distance=6mm] {node [leaf] {$e_{26}$}}
          child [sibling distance=10mm] {node [parallel] {$p$}
          	child [sibling distance=6mm] {node [series] {$s$}
          	 child [sibling distance=6mm] {node [leaf] {$e_{27}$}}
          	child [sibling distance=6mm] {node [leaf] {$e_{28}$}}}
            child [sibling distance=6mm] {node [leaf] {$e_{29}$}}
          }
        }
        child [sibling distance=22mm] {node [series] {$s$}
          child [sibling distance=10mm] {node [parallel] {$p$}
            child [sibling distance=6mm] {node [leaf] {$e_{30}$}}
            child [sibling distance=10mm] {node [series] {$s$}
              child [sibling distance=6mm] {node [leaf] {$e_{31}$}}
              child [sibling distance=6mm] {node [leaf] {$e_{32}$}}
            }
          }
          child [sibling distance=6mm] {node [leaf] {$e_{33}$}}
        }
      }
    }
    child [sibling distance=0mm] {node [leaf] {$e_{34}$}}
    child [sibling distance=15mm] {node [series] {$s$}
      child [sibling distance=6mm] {node [leaf] {$e_{35}$}}	
      child [sibling distance=10mm] {node [parallel] {$p$}
        child [sibling distance=6mm] {node [leaf] {$e_{36}$}}
        child [sibling distance=10mm] {node [series] {$s$}
          child [sibling distance=6mm] {node [leaf] {$e_{37}$}}
          child [sibling distance=6mm] {node [leaf] {$e_{38}$}}
        }
      }
    }
  }
;
\end{scope}
\end{tikzpicture}
		\caption{The series--parallel tree corresponding to the graph in Fig.~\ref{sp-graph}.}\label{sp-tree}
	\end{figure}
	
	In Section~\ref{lemmas} we will prove some useful lemmas, in Section~\ref{greater}
	we give a characterization of minimally $t$-tough sp-graphs with $\frac{1}{2}<
	t \leq 1$, in Section~\ref{half} we give a characterization of minimally $\frac{1}{2}$-tough
	sp-graphs, the main result is stated in Theorem~\ref{mainthm}.  Finally,  Section~\ref{open}  presents some open questions.

	\section{Useful Lemmas}\label{lemmas}
	
	The first lemma establishes that we can restrict our attention to simple graphs. While it might seem natural to assume that we are working exclusively with simple graphs, sp-graphs are not always simple. Although sp-graphs do not contain loops, they can include parallel edges.
	\begin{lem}
		\label{l5} If there are loops or parallel edges in graph $G$, then $G$ is
		not minimally $t$-tough for any value of $t>0$.
	\end{lem}
	\begin{proof}
		It is easy to see that removing a loop does not change the toughness.
		
		Let $e_{1}$ and $e_{2}$ be parallel edges, and let $S$ be a cutset. If $S$ contains an endpoint of $e_{1}$ (or $e_{2}$) then $c(G-e_{1}-S)=c(G-S)$, so
		$\tau(G-e_{1})=\tau(G)$, thus $G$ is not minimally $t$-tough.

		Moreover, if $S$ is not incident with $e_{1}$ and $e_{2}$ then
		$c(G-e_{1}-S)=c(G-S)$, since $e_{1}$ connects two vertices that are in the same
		component of $G-S$, i.e.~they are connected by $e_{2}$. Therefore, $G$ is not
		minimally $t$-tough.
	\end{proof}

	The next lemma provides a useful inequality involving the \emph{mediant} of two fractions.  
	The \emph{mediant} of \(\frac{a}{b}\) and \(\frac{c}{d}\) is defined as \(\frac{a+c}{b+d}\).  
	Although it is easy to prove Lemma~\ref{mediant}, it is interesting that it was already published in  1484 \cite{mediant}, and a straightforward ``proof without words'' can be found in \cite{mediantproof}.

	\begin{lem}
		[Mediant Inequality]\label{mediant} Let $a,b,c,d$ be positive integers. If $\frac{a}{b}
		\leq\frac{c}{d}$ then
		\[
		\frac{a}{b}\leq \frac{a+c}{b+d}\leq\frac{c}{d}.
		\]
		Moreover, equalities hold if and only if $\frac{a}{b}=\frac{c}{d}$.
	\end{lem}
	
	In the following, we show that if a graph has a cutset of size 1 or 2, then its
	tough sets satisfy some useful properties.
	
	\begin{lem}
		\label{l1} Let $G=G_{1}\cup G_{2}$, where $V(G_{1}) \cap V(G_{2})= \{v\}$. If
		$S$ is a tough set that does not contain $v$, then $S$ must be fully
		contained in either $V(G_{1})-\{v\}$ or $V(G_{2})-\{v\}$. (See Fig. \ref{fig:1}.)
	\end{lem}
	Note that $G$ is essentially the series join of $G_{1}$ and $G_{2}$, but we do
	not assume that $G_{1}$ and $G_{2}$ are sp-graphs.
	
	\begin{proof}
		Let $S$ be a tough set of $G$ so that $v\notin S$, and let
		$S_{1}=S\cap V(G_{1})$ and $S_{2}=S\cap V(G_{2})$. Suppose, to the contrary, that
		$S_{1}\not=\emptyset$ and $S_{2}\not=\emptyset$. \textit{ (We will use the notation
			$S_{1}$ and $S_{2}$ in a similar way in the next few lemmas, without giving
			the explicit definition.)}
		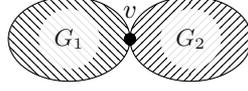
\begin{figure}[ht]
			\centering
			\begin{tikzpicture}[scale=0.8]
				\fill[pattern=north west lines] (-1,0) ellipse (1cm and 7mm); \fill[pattern=north
				east lines] (1,0) ellipse (1cm and 7mm);
				\begin{scope}
					\clip (-1,0) ellipse (1cm and 7mm); \clip (1,0) ellipse (1cm and 7mm);
					\fill[white] (-.5,-1) rectangle (.5,1);
				\end{scope}
				\draw (-1,0) ellipse (1cm and 7mm); \draw (1,0) ellipse (1cm and 7mm); \node[fill=white, circle,opacity=0.9]
				at (-1,0) {\small{\bf $G_{1}$}}; \node[fill=white,circle,opacity=0.9] at
				(1,0) {\small{\bf $G_{2}$}};
				\node[circle,fill=black,inner sep=1.5pt,draw] (a) at (0,0) {}; \node[yshift=10pt]
				at (a) {$v$};
			\end{tikzpicture}
			\caption{The structure of the graphs studied in Lemma~\ref{l1}.}
			\label{fig:1}
		\end{figure}
		Since $S$ is a tough set in $G$, we have $\tau(G)= \frac{|S|}{c(G-S)}$. We claim
		that $S_{i}$ is a cutset of $G$ i.e.,~$c(G-S_{i})\geq 2$ for $i=1,2$. If
		$S_{1}$ were not a cutset, then $c(G-S_{2})=c(G-S)$; thus
		\[
		\frac{|S_{2}|}{c(G-S_{2})}< \frac{|S|}{c(G-S)}=\tau(G),
		\]
		which would contradict the fact that $S$ is a tough set. By symmetry, both $S_{1}$ and $S_{2}$ are cutsets of $G$, which implies that
		$\tau(G)\le \frac{|S_{i}|}{c(G-S_{i})}$ for $i=1,2$. It is easy to see that $c
		(G-S)=c(G-S_{1})+c(G-S_{2}) -1>0$ holds. Without loss of generality, we can assume
		that $\frac{|S_{2}|}{c(G-S_{2})}\leq \frac{|S_{1}|}{c(G-S_{1})}$. Now the
		Mediant Inequality (Lemma~\ref{mediant}) implies that
		\[
		\frac{|S_{2}|}{c(G-S_{2})}\leq \frac{|S_{1}|+|S_{2}|}{c(G-S_{1})+c(G-S_{2})}
		\leq \frac{|S_{1}|}{c(G-S_{1})},
		\]
		and hence
		\[
		\tau(G)\le\frac{|S_{2}|}{c(G-S_{2})}< \frac{|S_{1}|+|S_{2}|}{c(G-S_{1})+c(G-S_{2})-1}
		=\tau(G),
		\]
		a contradiction.
	\end{proof}
	
	\begin{lem}
		\label{l2} Let $G = G_{1} \cup G_{2}\cup v_{1} v_{2}$, where
		$v_{i} \in V(G_{i})$ and $V(G_{1})\cap V(G_{2})=\emptyset$. If $S$ is a tough
		set, then $S$ must be fully contained in either $V(G_{1})$ or $V(G_{2})$.
	\end{lem}
	Note that, unlike in Lemma~\ref{l1}, now $S$ may contain $v_{1}$ or $v_{2}$ or
	even both.
	\begin{figure}[ht]
		\centering
		\begin{tikzpicture}[scale=0.8]
			\fill[pattern=north west lines] (-1,0) ellipse (1cm and 7mm); \fill[pattern=north
			east lines] (2,0) ellipse (1cm and 7mm);
			\begin{scope}
				\clip (-1,0) ellipse (1cm and 7mm); \clip (2,0) ellipse (1cm and 7mm);; \fill[white]
				(-.5,-1) rectangle (.5,1);
			\end{scope}
			\draw (-1,0) ellipse (1cm and 7mm); \draw (2,0) ellipse (1cm and 7mm); \node[fill=white,circle,
			opacity=0.9] at (-.9,0) {\small{\bf $G_{1}$}}; \node[fill=white,circle,
			opacity=0.9] at (2,0) {\small{\bf $G_{2}$}}; \node[circle,fill=black,inner
			sep=1.5pt,draw] (a) at (0.008,0.01) {}; \node[yshift=7pt,xshift=5pt] at (a)
			{$v_{1}$}; \node[circle,fill=black,inner sep=1.5pt,draw] (b) at (1,-0.009)
			{}; \node[yshift=7pt,xshift=-5pt] at (b) {$v_{2}$}; \node (v_1) at (0.008,0.01)
			{}; \node (v_2) at (1,-0.009) {}; \draw[] (a) -- (b);
		\end{tikzpicture}
		\caption{The structure of the graphs studied in Lemma~\ref{l2}.}
		\label{fig:2}
	\end{figure}
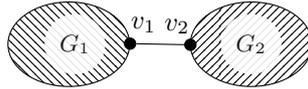
	\begin{proof}
		If $v_{1}\notin S$ or $v_{2}\notin S$, then we can apply Lemma~\ref*{l1} by setting
		$v_{1}$ or $v_{2}$ as $v$. Hence, we can assume that $v_{1},v_{2}\in S$.
		
		Suppose to the contrary that $S_{1}=S \cap V(G_{1})\neq \emptyset$ and $S_{2}
		=S \cap V(G_{2})\neq \emptyset$. Since $S$ is a tough set in $G$, we have $\tau
		(G)= \frac{|S|}{c(G-S)}$. Since $S$ is a tough set, it is easy to see that $V
		(G_{1})\nsubseteq S$ and $V(G_{2})\nsubseteq S$. Therefore $S_{1}$ and
		$S_{2}$ are both cutsets in $G$, which implies that $\tau(G)\le \frac{|S_{i}|}{c(G-S_{i})}$
		holds for $i=1,2$.
		
		Note that in this case $c(G-S)= c(G-S_{1})+c(G-S_{2}) -2$. Without loss of
		generality, we can assume that
		$\frac{|S_{2}|}{c(G-S_{2})}\leq \frac{|S_{1}|}{c(G-S_{1})}$. Now the Mediant
		Inequality (Lemma~\ref{mediant}) implies that
		\[
		\frac{|S_{2}|}{c(G-S_{2})}\leq \frac{|S_{1}|+|S_{2}|}{c(G-S_{1})+c(G-S_{2})}
		\leq \frac{|S_{1}|}{c(G-S_{1})}
		\]
		and hence
		\[
		\tau(G)\le\frac{|S_{2}|}{c(G-S_{2})}< \frac{|S_{1}|+|S_{2}|}{c(G-S_{1})+c(G-S_{2})-2}
		=\tau(G),
		\]
		a contradiction.
	\end{proof}
	
	\begin{lem}
		\label{l3} Let $G = G_{1} \cup G_{2}$, where
		$V(G_{1}) \cap V(G_{2})=\{u,v\}$. Assume that $\tau(G) < 1$. If $S$ is tough
		set of minimum size, and $S\cap\{u,v\}=\emptyset$, then $S$ must be fully
		contained in either $V(G_{1})-\{u,v\}$ or $V(G_{2})-\{u,v\}$. Moreover, if
		$S$ is a tough set, $S\cap\{u,v\}=\emptyset$, and $|S_{1}|,|S_{2}|>0$ (hence
		$S$ is not minimum size by the first claim), then $\frac{|S_{1}|}{c(G-S_{1})}
		= \frac{|S_{2}|}{c(G-S_{2})}=\tau(G)$.
	\end{lem}
	Note that there may be an edge in $G$ between $u$ and $v$. In this case, we
	can assume that this edge belongs to $G_{1}$.
	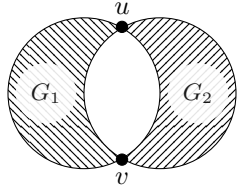
\begin{figure}[ht]
		\centering
		\begin{tikzpicture}[scale=1]
			\fill[pattern=north west lines] (-0.5,0) circle (1cm); \fill[pattern=north
			east lines] (0.5,0) circle (1cm);
			\begin{scope}
				\clip (-0.5,0) circle (1cm); \clip (0.5,0) circle (1cm); \fill[white] (-.5,-1)
				rectangle (.5,1);
			\end{scope}
			\draw (-0.5,0) circle (1cm); \draw (0.5,0) circle (1cm); \node[fill=white,circle,
			opacity=0.9] at (-1,0) {\small{\bf $G_{1}$}}; \node[fill=white,circle,
			opacity=0.9] at (1,0) {\small{\bf $G_{2}$}}; \node[circle,fill=black,inner
			sep=1.5pt,draw] (a) at (0,.88) {}; \node[yshift=7pt] at (a) {$u$}; \node[circle,fill=black,inner
			sep=1.5pt,draw] (b) at (0,-.88) {}; \node[yshift=-7pt] at (b) {$v$};
		\end{tikzpicture}
		\caption{The structure of the graphs studied in Lemma~\ref{l3}.}
		\label{fig:3}
	\end{figure}
	\begin{proof}
		Suppose to the contrary that $S$ is a tough set in $G$, such that
		$S\cap\{u,v\}=\emptyset$, but $S_{1}=S \cap V(G_{1})\neq \emptyset$ and $S_{2}
		=S \cap V(G_{2})\neq \emptyset$. This implies that $|S|\ge 2$. Since $S$ is a
		tough set in $G$, we have $\tau(G)= \frac{|S|}{c(G-S)}$.
		
		If $c(G-S_{1})=1$ and $c(G-S_{2})=1$, then each vertex of $G-S_{1}$ and $G-S_{2}$
		is connected with a path to either $u$ or $v$. This implies that
		$c( G-S)\leq 2$, hence $\tau(G)\geq 1$, a contradiction.
		
		If $c(G-S_{1})=1$ and $c(G-S_{2})\geq 2$, with $|S_{1}|, |S_{2}|\geq 1$,
		then $c( G-S)\leq c(G-S_{2})+1$ holds, since in $G-S_{2}$ the vertices $u$ and
		$v$ belong to the same component, but in $G-S_{1}$ each vertex is connected
		to either $u$ or $v$. So if we remove $S_{1}$ from $G-S_{2}$, then the
		number of components may only increase by $1$.
		
		If $\frac{|S_{2}|}{c(G-S_{2})}<1$, then the Mediant Inequality (Lemma
		\ref{mediant}) implies that
		\[
		\frac{|S_{2}|}{c(G-S_{2})}< \frac{|S_{2}|+1}{c(G-S_{2})+1}<\frac{1}{1}.
		\]
		Therefore we have
		\[
		\begin{split}
			\frac{|S_{2}|}{c(G-S_{2})}< \frac{|S_{2}|+1}{c(G-S_{2})+1}\leq \frac{|S_{1}|+|S_{2}|}{c(G-S_{2})+1}
			\\ \leq\frac{|S|}{c( G-S)}= \tau(G),
		\end{split}
		\]
		so $S$ cannot be a tough set.\\
		
		If $\frac{|S_{2}|}{c(G-S_{2})}\ge 1$, then the Mediant Inequality (Lemma
		\ref{mediant}) implies that
		\[
		\frac{1}{1}\le \frac{|S_{2}|+1}{c(G-S_{2})+1}\le\frac{|S_{2}|}{c(G-S_{2})}.
		\]
		Now we have
		\[
		\begin{split}
			1 \le \frac{|S_{2}|+1}{c(G-S_{2})+1}\leq \frac{|S_{1}|+|S_{2}|}{c(G-S_{2})+1}
			\leq\\ \frac{|S|}{c( G-S)}= \tau(G),
		\end{split}
		\]
		contradicting the assumption.
		
		Hence, we can assume that $c(G-S_{i})\geq 2$ holds for $i=1,2$. Since both $S
		_{1}$ and $S_{2}$ are cutsets in $G$, we have
		$\tau(G)\le \frac{|S_{i}|}{c(G-S_{i})}$ for $i=1,2$.
		
		Consider the components of $G-S$.
		\begin{enumerate}
			\item[(i)] The components that are disjoint from $V(G_{1})$, i.e.~components
			of $G-S_{2}$ not containing $u$ or $v$. There are $c(G-S_{2})-1$ such
			components.
			
			\item[(ii)] The components that are disjoint from $V(G_{2})$, i.e.~components
			of $G-S_{1}$ not containing $u$ or $v$. There are $c(G-S_{1})-1$ such
			components.
			
			\item[(iii)] \label{type3} The components of $u$ and $v$. These may be
			separate components, but if there is a path in $G-S$ between $u$ and $v$,
			then they are the same. So there are $1$ or $2$ such components.
		\end{enumerate}
		Thus
		\[
		\begin{split}
			c( G-S)=c(G-S_{1})-1 +c(G-S_{2})-1 +1=\\ c(G-S_{1})+c(G-S_{2})-1,
		\end{split}
		\]
		if $u$ and $v$ are in one component of $G-S$, and
		\[
		\begin{split}
			c( G-S)=c(G-S_{1})-1 +c(G-S_{2})-1 +2=\\ c(G-S_{1})+c(G-S_{2}),
		\end{split}
		\]
		if $u$ and $v$ are in different components of $G-S$. Therefore we have
		\[
		\begin{split}
			c(G-S_{1})+c(G-S_{2})\geq c( G-S)\geq\\ c(G-S_{1})+ c(G-S_{2})-1.
		\end{split}
		\]
		
		Without loss of generality, we can assume that
		$\frac{|S_{1}|}{c(G-S_{1})}\leq \frac{|S_{2}|}{c(G-S_{2})}$. Now the Mediant
		Inequality (Lemma~\ref{mediant}) implies that
		\begin{equation}
			\label{eq:2.1}\frac{|S_{1}|}{c(G-S_{1})}\leq \frac{|S_{1}|+|S_{2}|}{c(G-S_{1})+
				c(G-S_{2})}\leq\frac{|S_{2}|}{c(G-S_{2})},
		\end{equation}
		and hence
		\begin{equation}
			\label{eq:2.2}\frac{|S_{1}|}{c(G-S_{1})}\leq \frac{|S|}{c( G-S)},
		\end{equation}
		so $S$ cannot be a minimum-size tough set. If the inequality in \eqref{eq:2.2}
		is strict, then $S$ cannot even be a tough set. If equality holds in
		\eqref{eq:2.2}, then there is also equality in \eqref{eq:2.1}, so
		\[
		\frac{|S_{1}|}{c(G-S_{1})}=\tau(G)= \frac{|S_{2}|}{c(G-S_{2})}
		\]
		holds, completing the proof.
	\end{proof}
	\begin{lem}
		\label{l4} Let $G=G_{1}\cup G_{2}\cup u_{1}u_{2}\cup v_{1}v_{2}$ where $\{u_{i}
		,v_{i}\} \subseteq V(G_{i})$ holds for $i=1,2$. Assume that $\tau(G) < 1$.
			If $S$ is a minimum-size tough set, then $S$ must be fully contained in either
		$V(G_{1})$ or $V(G_{2})$ (see Fig.~\ref{fig:4}). Moreover, if $S$ is a tough
		set with $|S_{1}|,|S_{2}|>0$ (hence $S$ is not minimal), then $\frac{|S_{1}|}{c(G-S_{1})}
		= \frac{|S_{2}|}{c(G-S_{2})}=\tau(G)$.
	\end{lem}
	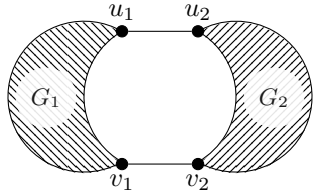
\begin{figure}[ht]
		\centering
		\begin{tikzpicture}[scale=1]
			\fill[pattern=north west lines] (-0.5,0) circle (1cm); \fill[pattern=north
			east lines] (1.5,0) circle (1cm);
			\begin{scope}
				\clip (0.5,0) circle (1cm); \fill[white] (-.5,-1) rectangle (1.5,1);
			\end{scope}
			\draw (0,.88) arc (60:300:1) arc (240:120:1);
			\draw (1,.88) arc (60:-60:1) arc (240:480:1); \node[fill=white,circle,
			opacity=0.9] at (-1,0) {\small{\bf $G_{1}$}}; \node[fill=white,circle,
			opacity=0.9] at (2,0) {\small{\bf $G_{2}$}}; \node[circle,fill=black,inner
			sep=1.5pt,draw] (a1) at (0,.88) {}; \node[yshift=7pt] at (a1) {$u_{1}$}; \node[circle,fill=black,inner
			sep=1.5pt,draw] (b1) at (0,-.88) {}; \node[yshift=-7pt] at (b1) {$v_{1}$};
			\node[circle,fill=black,inner sep=1.5pt,draw] (a2) at (1,.88) {}; \node[yshift=7pt]
			at (a2) {$u_{2}$}; \node[circle,fill=black,inner sep=1.5pt,draw] (b2) at (1,-.88)
			{}; \node[yshift=-7pt] at (b2) {$v_{2}$}; \draw (a1)--(a2); \draw (b1)--(b2);
		\end{tikzpicture}
		\caption{The structure of the graphs studied in Lemma~\ref{l4}.}
		\label{fig:4}
	\end{figure}
	\begin{proof}
		We can argue similarly to that in the proof of Lemma~\ref{l3}. Notice that
		in $G-S$, the vertices $u_{1}$ and $u_{2}$ (resp. $v_{1}$ and $v_{2}$)
		cannot belong to different components. 
        On the other hand, one or both of them may not be in the graph $G-S$. 
        So there may be 0 or 1 component containing
		$u_{1}$ or $u_{2}$ (resp. $v_{1}$ or $v_{2}$). Therefore, there may be
		$0, 1$, or $2$ components of type (iii).
		
		This leads to the inequality
		\[
		\begin{split}
			c(G-S_{1})+c(G-S_{2})\geq c( G-S) \geq\\ c(G-S_{1})+c(G-S_{2})-2.
		\end{split}
		\]
		
		Now, the proof can be finished in the same way.
	\end{proof}
	
	We call a path $u,w_{1},w_{2},\ldots, w_{k},v$ a \emph{thread} if $k\ge 2$ and the degree of all $w_i$ vertices is 2.  A path $u,w,v$ is called a \emph{short-thread} if the degree  of $w$ is 2. Our next goal is to show that it is enough to consider graphs that do not	contain any threads (only possibly some short-threads).
	
	\begin{claim}\label{claim} 
		Let \( G \) be a minimally \( t \)-tough graph. For any edge \( e = uv \) in \( G \), no tough set of \( G - e \) can contain either \( u \) or \( v \).
	\end{claim}
	
	\begin{proof}
		Let $S$ be a tough set of $G-e$. If $u\in S$ or $v\in S$, then $(G-e)-S$ is the same graph as $G-S$, so removing $e$ does not change the toughness, a contradiction.
	\end{proof}
	
	\begin{lem}
		\label{p1} Let $t<1$, and $u,w_{1},w_{2},\ldots, w_{k},v$ be a thread or a short-thread in
		a   graph $G$ with toughness $t$ such that  $G-\{w_{1},w_{2},\ldots, w_{k}\}$ is connected. If $S$
		is a minimum-size tough set of $G$, then
		$S \cap \{w_{1},w_{2},\ldots, w_{k}\}=\emptyset.$
	\end{lem}


	\begin{proof}
		First, to prove the claim for the threads, we can apply Lemma~\ref{l4} by setting $u_{1}u_{2}=uw_{1}$ and $v_{1}v_{2}=w_{k}
		v$. Let $G_{1}$ be the path $w_{1},w_{2},\ldots, w_{k}$, and $G_{2}$ the rest of $G$. If $\emptyset\not=S\subseteq V(G_{1})$ is a minimum-size tough set, then it
		is easy to see that $\frac{|S_{1}|}{c(G-S_{1})}\ge 1$, which contradicts the assumption
		$t<1$. Hence $S$ is entirely contained in $G_2$ proving the claim.
		
		Now we prove the claim for the short-threads. We can apply Lemma~\ref{l3} by setting $G_1$ to be the path $u,w,v$ and $G_2$ the rest of the graph. If $\emptyset\not=S$ is a minimum-size tough set that contains neither \(u\) nor \(v\), then clearly it cannot contain $w$. If it contains $w$ and at least one of $u$ or $v$ then restoring $w$ from $S$ cannot decrease the number of components, so $\frac{|S-w|}{c(G-(S-w))}\le\frac{|S|}{c(G-S)}$, a contradiction.
              \end{proof}
  
    \begin{lem}
		\label{p2} Let $t<1$, and $u,w_{1},w_{2},\ldots, w_{k},v$ be a thread or a short-thread in
		a   graph $G$ with toughness $t$ such that $u$ is a cut-vertex and $G-\{w_{1},w_{2},\ldots, w_{k}\}$ is disconnected. Then there exists
		a tough set $S$ with $S \cap \{w_{1},w_{2},\ldots, w_{k}\}=\emptyset$.
     \end{lem}
              
\begin{proof}
	Since $u$ is a cut-vertex, we have $\tau(G)\leq \frac{1}{2}$.
 If	$\tau(G)= \frac{1}{2}$ then $\{u\}$ is the desired tough set.
	
	Now suppose that $\tau(G)<\frac{1}{2}$, and let $S$ be a tough
	set of $G$. We show that no vertex of degree at most two belongs
	to $S$.
	
	Suppose, to the contrary, that $x\in S$ and $d_G(x)\leq 2$.
	Since $S$ is a tough set,
	$
	\frac{|S|}{c(G-S)}=\tau(G)<\frac{1}{2},
	$
	and consequently $c(G-S)>2|S|$.
	
	Since $G$ is connected and $d_G(x)\leq 2$, the case
	$|S|=1$ would imply $c(G-S)\leq 2$, contradicting
	$c(G-S)>2|S|$. Hence $|S|\geq 2$. Consider $S'=S-x$. Since $d_G(x)\leq 2$, restoring $x$
	can merge at most two components of $G-S$. Therefore
	\[
	c(G-S')\geq c(G-S)-1>2|S|-1\geq 1.
	\]
	In particular, $S'$ is still a cutset. Moreover, using  $c(G-S)>2|S|>|S|$, we obtain
	\[
	\frac{|S'|}{c(G-S')}
	\leq \frac{|S|-1}{c(G-S)-1}
	<\frac{|S|}{c(G-S)}=\tau(G),
	\]
	a contradiction.
	
	Hence no vertex of degree at most two belongs to $S$.
	In particular, since $w_1,w_2,\ldots,w_k$ all have degree two,
	we obtain
	$
	S\cap\{w_1,w_2,\ldots,w_k\}=\emptyset$,
		as required.
\end{proof}
  
\begin{lem}
	\label{p3} Suppose that $G'$ is obtained from $G$ with toughness $t<1$ by replacing a thread with a short-thread.
	\begin{itemize}
		\item[(i)]\label{cl1} 
		$\tau(G)= \tau(G')$.
 		\item[(ii)] $G'$ is
		minimally $t$-tough if and only if $G$ is minimally $t$-tough.
	\end{itemize}
    \end{lem}
              
	\begin{proof} Suppose that we replace the thread
		$u,w_1,w_2,\ldots,w_k,v$ with the short-thread $u,w,v$.
		The first case to consider is when $G$ is itself a path.
		In this case, both claims hold, since every path of
		length at least 3 is minimally $\frac12$-tough.
	
	So we can assume that $u$ or $v$ has at least one more neighbor outside of the thread. This implies that either $u$ or $v$  is a cut-vertex or $G-\{w_{1},w_{2},\ldots, w_{k}\}$ is connected.
		
		\begin{proof}[Proof of (i).]
			By Lemma~\ref{p1} and \ref{p2} there exists a tough set $S$ in $G$ with $S \cap
			\{w_{1},w_{2},\ldots ,w_{k}\}=\emptyset$. Thus $c(G'-S)=c(G-S)$, since in $G
			-S$ all $w_{i}$ belong to the same component, which corresponds to the component
			of $G'-S$ that contains $w$. This implies that
			$\tau(G)=\frac{|S|}{c(G-S)}=\frac{|S|}{c(G'-S)}\geq \tau(G')$.
			
			For the other direction, we observe that \(\tau(G')\leq\tau(G)=t<1\). If \(G'-w\) is connected, then Lemma~\ref{p1} guarantees the existence of a tough set of \(G'\) avoiding \(w\). Otherwise, at least one of \(u,v\) is a cut-vertex, and after interchanging \(u\) and \(v\), if necessary, Lemma~\ref{p2} applies. Thus, \(G'\) contains a tough set \(S'\) with \(w\notin S'\). Since \(c(G-S')=c(G'-S')\), we obtain \(\tau(G')\geq\tau(G)\), completing the proof of the claim.			
                      \end{proof}
                      
                      
		\begin{proof}[Proof of (ii).]
			Let $e$ be an arbitrary edge of $G$. If $e$ is not on the thread, then we can apply part~(i) for $G-e$ as well.
			This implies that $\tau(G-e)= \tau(G'-e)$, so deleting $e$ decreases
			the toughness of $G$ iff it decreases the toughness of $G'$.
			

If $e$ is an edge of the thread, we distinguish two
directions.

Suppose first that $G$ is minimally $t$-tough.
Consider $H=G-uw_1$. By minimality, $\tau(H)<t$.
If $H$ is disconnected, then $G'-uw$ is also disconnected.
Otherwise, $v$ is a cut-vertex of $H$, and hence
$\tau(H)\leq 1/2$.

If $\tau(H)=1/2$, then $\{v\}$ is a tough set of $H$.
If $\tau(H)<1/2$, then the path
$v,w_k,w_{k-1},\ldots,w_2,w_1$ is a thread in $H$.
Since $v$ is a cut-vertex of $H$ and
$H-\{w_2,\ldots,w_k\}$ is disconnected, we can apply
Lemma~\ref{p2}, obtaining a tough set $S$ of $H$
avoiding $w_2,\ldots,w_k$.
Moreover, $w_1$ cannot belong to $S$: since $w_1$
has degree one, restoring it after deleting $S$ cannot
decrease the number of components. Thus removing $w_1$
from $S$ would give a cutset with a smaller toughness
ratio, a contradiction.

Since $H$ and $G'-uw$ differ only in the length of a
pendant path attached at $v$, we have
$
c(H-S)=c((G'-uw)-S).
$
Consequently,
$
\tau(G'-uw)\leq\tau(H)<t.
$
The same argument, with the roles of $u$ and $v$
interchanged, shows that $\tau(G'-wv)<t$.

In the direction where $G'$ is assumed to be minimally
$t$-tough, we apply the proof of (i) starting from $G'-e$.
Since $\tau(G'-e)<t<1$, the same argument gives
$\tau(G-e)=\tau(G'-e)$.
Let $e=w_iw_{i+1}$ be an internal edge of the original
thread. Since $G'$ is minimally $t$-tough,
$\tau(G'-uw)<t$.

If $G'-uw$ is disconnected, then $G-e$ is also
disconnected. Otherwise, let $S$ be a tough set of
$G'-uw$. By Claim~\ref{claim}, neither $u$ nor $w$
belongs to $S$.

In $G-e$, the original thread splits into two paths,
one attached at $u$ and the other at $v$. Since $u\notin S$,
the path attached at $u$ does not create an additional
component. The path attached at $v$ has the same effect
on the number of components as the pendant path $wv$ in
$G'-uw$. Therefore
$c((G-e)-S)=c((G'-uw)-S)$.

This implies that
\[
\tau(G-e)\leq
\frac{|S|}{c((G-e)-S)}
=\tau(G'-uw)<t.
\]

For $e=uw_1$, put $H=G-uw_1$ and $H'=G'-uw$. If $H'$ is disconnected, then $H$ is also disconnected,
and hence $\tau(H)=0<t$.
Otherwise, by Claim~\ref{claim}, a tough set $S$ of $H'$ avoids
both $u$ and $w$. Since $H$ and $H'$ differ only in
the length of the pendant path attached at $v$, we have
$c(H-S)=c(H'-S)$. Consequently,
\[
\tau(G-uw_1)\leq
\frac{|S|}{c(H-S)}
=\tau(H')<t.
\]
The other end edge is handled symmetrically.
Thus, deleting any edge of $G$ decreases its toughness.
This completes the proof.
		\end{proof}
	\end{proof}

 	\begin{defn}
		Let $r(G)$, \emph{the reduced graph of $G$}, be the graph obtained from $G$ by
		repeatedly replacing a thread with a path of
		length $2$ until there is no thread remains in the graph.
	\end{defn}
	
	Lemma~\ref{p3} proves that if $t<1$, then $G$ is minimally $t$-tough if and only
	if $r(G)$ is minimally $t$-tough. Therefore, if we can characterize the reduced minimally tough graphs, then we can also characterize all such graphs.

	The edges $e$, whose parent in the sp-tree (i.e.~its first join) is a parallel join, are special.
	Note that in this case, the siblings of $e$ cannot be edges since we assume that
	there are no parallel edges. So, all siblings of $e$ must be series nodes. We will
	see that some edges have different properties from others, so we define two
	special types of edges. {The criteria below are defined using the sp-tree in Fig.~\ref{sp-tree}.}
	
	\begin{defn}
		An edge $e$ is a {\em leap-edge} if its parent is a parallel node, which is the
		root of the sp-tree, and it has only one sibling (which is a series node).
		An edge is a {\em jump-edge} if its parent is a parallel node that is not the
		root or the parent is the root and has at least two other children, all of which are series nodes
		(i.e.~it is not a leap-edge).
	\end{defn}
	
Notice that if two edges are parallel, then both of them are either leap-edges or jump-edges.  { Thus, to classify any labeled edge $e_i$ in Fig.~\ref{sp-graph}, locate its leaf in Fig.~\ref{sp-tree} and inspect the parent and sibling nodes as specified above. The edges $e_{15}, e_{29}, e_{30}, e_{34}, e_{36}$ are all jump-edges. There is no leap-edge in this example, $e_{34}$ is a jump-edge. However,
when the subgraph formed by \(e_{34},\ldots,e_{38}\) in Fig.~\ref{sp-graph} is considered as a series--parallel graph on its own, then \(e_{34}\) is a leap-edge in this graph.}

      \begin{lem}
          	\label{j1}
          	Let $G$ be a graph with $\frac{1}{2} \le \tau(G) \le 1$, and let $e$ be a jump-edge in $G$. Then $\tau(G - e) = \tau(G)$ in either of the following cases:
          	\begin{itemize}
          		\item[(i)] $\tau(G) > \frac{1}{2}$; or
          		\item[(ii)] $\tau(G) = \frac{1}{2}$ and the grandparent of the leaf corresponding to $e$ in the sp-tree (a series node) is not the root.
          	\end{itemize}
      \end{lem}

	\begin{proof}
		Let $e=s_{1}t_{1}$; thus the parent of $e$ is a parallel node with joining
		vertices $s_{1}$ and $t_{1}$. 
		
		\subsubsection*{Case 1: $e$ has a grandparent:}
		The grandparent must be a series node. If it is the root, then removing a
		joining vertex gives at least $2$ components, thus $\tau(G) \leq \frac{1}{2}$,
		implying that neither (i) nor (ii) holds. Therefore, we can assume that $e$
		has a great-grandparent that is a parallel node.
		
		Suppose that $\tau(G-e)< \tau(G)$ and let $S$ be a minimum-size tough set for $G-e$. If $s
		_{1} \in S$ or $t_{1} \in S$, then $C((G-e)-S)= C(G-S)$ implying $\tau (G-e) \geq
		\tau(G)$. Otherwise, the only difference between the components of $(G-e)-S$ and
		$G-S$ is that in $(G-e)-S$, the vertices $s_{1}$ and $t_{1}$ may be in
		different components, but $e$ connects them in $G-S$.
		
		$\{s_{1},t_{1}\}$ is clearly a cutset. Let $G_{1}$ denote the subgraph containing
		the series node sibling of $e$, and $G_{2}$ the rest of the graph. Applying
		Lemma~\ref{l3} by setting $u=s_{1}$ and $v=t_{1}$ gives that $S$ is contained in either
		$(G_{1}-s_{1}-t_{1})$ or $(G_{2}-s_{1}-t_{1})$.
		
		If $S$ is contained in $(G_{1}-s_{1}-t_{1})$, then there is a path in $G-e$
		between $s_{1}$ and $t_{1}$ through $G_{2}$, thus $C(G-e)=C(G)$, which
		implies $\tau (G-e) \geq \tau(G)$, a contradiction.
		
		If $S$ is contained in $G_{2} -s_{1}-t_{1}$, then there is a path in $G-e$ between
		$s_{1}$ and $t_{1}$ through $G_{1}$, thus $C(G-e)=C(G)$ which implies
		$\tau (G-e) \geq \tau(G)$, a contradiction again.
		
		\subsubsection*{Case 2: $e$ has no grandparent:}
		If the parent of $e$ is a parallel join, which is the root, but it has $\geq 2$
		siblings that are series nodes, then the argument is similar; just let
		$G_{1}, G_{2}$ be the subgraphs corresponding to the two series siblings.
	\end{proof}

			\section{Characterization of minimally $t$-tough sp-graphs with
			\texorpdfstring{$\frac{1}{2}< t \leq 1$}{1/2 < t <= 1}}\label{greater}
		
	\begin{thm}
		If $G$ is a minimally $1$-tough sp-graph, then $G$ is a cycle.
	\end{thm}
	\begin{proof}
		If the root of the sp-tree is a series node, then there is a cut-vertex, so $\tau(G)\leq \frac{1}{2}
		< 1$, a contradiction. So, the root must be a parallel node. It cannot have two children
		that are edges since they would be parallel.
		
		\subsection*{Case 1: The root has two children, an edge and a series node.}
		
		\subsubsection*{Subcase 1a: all children of the series node are edges.}
		In other words, the graph is a parallel join of an edge and a series join of some edges (i.e., a path) and an edge. So it is a cycle.
		\subsubsection*{Subcase 1b: the series node has a parallel node child.}
		Let us denote  this parallel node by $x$, and the two terminals of
		this join by $s_{x},t_{x}$. First, note that $x$ must have at least two children.
		If it has at least two series nodes, then $G-s_{x}-t_{x}$ has at least 3 components:
		two components from the series node children of $x$ and one more containing
		the rest of the vertices. This implies that $\tau(G)\le\frac{2}{3}<1$, a
		contradiction.
		
		$x$ cannot have more than one child, which is an edge, so it remains to be
		handled when $x$ has two children: an edge and a series node. However, this means
		that the edge is a jump-edge. Now Lemma~\ref{j1} implies that $G$ is not
		minimally $1$-tough.
		
		\subsection*{Case 2: the root has two children that are series nodes.}
		First, note that the root cannot have children that are edges, since they
		would be jump-edges, and then Lemma~\ref{j1} implies that $G$ cannot be minimally
		$1$-tough.
		
		\subsubsection*{Subcase 2a: all children of both series nodes are edges.}
		In other words, the graph is a parallel join of two series joins of some
		edges (i.e.~two paths). So it is a cycle.
		
		\subsubsection*{Subcase 2b: one of the series nodes has a parallel node child.}
		This case is essentially the same as Subcase 1b.
		
		\subsection*{Case 3: the root has at least 3 children that are series nodes.}
		Removing the two terminals of the root gives at least 3 components. This
		implies that $\tau(G)\le\frac{2}{3}<1$, a contradiction.
	\end{proof}
	
	Now, we are ready to state the main result of this section.
	\begin{thm}
		Let $G$ be an sp-graph with
		$\frac12<\tau(G)<1$. Then $G$ is minimally
		$\tau(G)$-tough if and only if it has no jump-edges.
	\end{thm}
	
\begin{proof}
	If there is a jump-edge, then the graph is not minimally tough by Lemma
	\ref{j1}. 
	
	Now, suppose there are no jump-edges. We show that for any edge $e$, we have $\tau(G-e)<\tau(G)$. In fact, we show
	that for any edge $\tau(G-e)\leq \frac{1}{2}< \tau(G).$
	
	\subsection*{Case 1: The parent of $e$ is a series node.}
	
	\subsubsection*{Subcase 1a: The parent of $e$ is the root.}
	{(For example $e_1$ in Fig.~\ref{sp-tree}.)} In this case, $G-e$ has $2$ components. This is true even in the special cases
	when $G$ is just an edge or if $e$ is incident to one of the terminals of the
	series node. In these cases, one or both of the components is a single vertex.
	This implies that $\tau(G-e)=0<\frac{1}{2}<\tau(G)$, a contradiction.
	\subsubsection*{Subcase 1b: The parent of $e$ is not the root.}
		{(For example $e_2$ in Fig.~\ref{sp-tree}.)} Therefore, the grandparent of $e$ must be a parallel node. Let $s$ and $t$ be the
	terminals of this parallel join. Also, let $G'$ be the subgraph of the series
	node with terminals $s$ and $t$. It is clear that all paths connecting $s$
	and $t$ that use the edges of $G'$ must contain $e$. If $e$ is not incident to
	either $s$ or $t$, then it is easy to see that $G'-e$ has two components, one
	containing $s$, the other containing $t$. It is clear that both $s$ and $t$ are cut
	vertices in $G-e$. On the other hand, if $e$ is incident to say $s$ (it
	cannot be incident to both), then $t$ is still a cut-vertex. This implies that
$\tau(G-e)\le\frac{1}{2}<\tau(G)$, a contradiction.
	
	\subsection*{Case 2: The parent of $e$ is a parallel node.}
		{(For example $e_{15}$ in Fig.~\ref{sp-tree}.)} If the parent is not the root, then it is a jump-edge, contradicting our
	assumption. If the parent is the root, then it may only have one sibling, 
    since otherwise, it would be a jump-edge. Thus, $G-e$ is a series join of at
	least two sp-graphs, so any of the internal join vertices are cut-vertices.
	This implies that $\tau (G-e) \leq \frac{1}{2}< \tau(G)$.
\end{proof}
	
	  \section{Characterization of minimally \texorpdfstring{$\frac{1}{2}$}{1/2}-tough sp-graphs}\label{half}
	
	In this section, we first characterize all reduced minimally $\frac{1}{2}$-tough sp-graphs. By Lemma~\ref{p3}, this will provide a characterization of all minimally $\frac{1}{2}$-tough sp-graphs. 
	
	For this, we define a few specific \textit{substructures} of the sp-tree. A substructure is like a subtree, but (only) the root may have other siblings in the sp-tree. A subtree is always a substructure as well. A special case is when the substructure is the whole graph.

\begin{figure}[h]
		\begin{center}
%
%
%
	\begin{tikzpicture}[
	every edge/.style = {draw,very thick},
	vertex/.style args = {#1 #2}{circle, 
		draw, 
		thick,
		minimum size=2pt,inner sep=-2pt,
		fill=black,
		label=#1:#2},
			leaf/.style={rectangle, draw=black, rounded corners=1mm, 
			text centered, anchor=north, text=black,text width=1em,inner sep=1.5pt},
		parallel/.style={circle, draw=black, fill=gray!30,
			text centered, anchor=north, text=black,text width=1em,inner sep=1pt},
		non/.style={circle, 
				text centered, anchor=north, text=black,text width=1em,inner sep=1pt},
		series/.style={circle, draw=black, fill=black!70, 
			text centered, anchor=north, text=white, text width=1em,inner sep=1pt},
		level distance=0.25cm, growth parent anchor=south,level/.style={sibling distance=10mm}]
		
	\begin{scope}
			\begin{scope}
			\node  at (-1,-.3) {$P_2$};
				\node (R) [series] {$s$}
				child[sibling distance=7mm] {node [leaf] {$e$}}
				child[sibling distance=7mm] {node [leaf] {$e$}}
			;
			\draw[dashed] (R) -- (1,-.8);
			\draw[dashed] (R) -- (-1,-.8);
		\end{scope}
		
		\begin{scope}[shift=({0,-2})]
	\path   
	node (A) [vertex=above $s_1$]     at ( -1, 0) {} 
	node (B) [vertex=above $v_1$]     at ( 0, 0) {} 
	node (C) [vertex=above $t_1$]     at ( 1, 0) {} 
	node (D)     at ( -1.7, 0) {} 
	node (E)      at ( 1.7, 0) {} 

	(A) edge (B)
	(B) edge (C)
	(D) edge[dashed] (A)
	(C) edge[dashed] (E)
	
;

	\end{scope}
	\end{scope}
	
		\begin{scope}[shift=({4,1})]
		\begin{scope}
			\node at (-1,-0.3) {$Q_2$};
			\node (R)[parallel] {$p$}
			child[sibling distance=15mm] {node [series] {$s$}
				child[sibling distance=7mm] {node [leaf] {$e$}}
				child[sibling distance=7mm] {node [leaf] {$e$}}}
			child[sibling distance=15mm] {node [leaf] {$e$}
				}
			;
			\draw[dashed] (R) -- (1.3,-.8);
			
		\end{scope}
		
		\begin{scope}[shift=({0,-3})]
			\path   
			node (A) [vertex=above $s_1$]     at ( -1, 0) {} 
			node (B) [vertex=above $v_1$]     at ( 0, .5) {} 
			node (C) [vertex=above $t_1$]     at ( 1, 0) {} 
		node (D)     at ( -0.7, -0.5) {} 
		node (E)      at ( 0.7, -0.5) {} 
			
			(A) edge (B)
			(B) edge (C)
			(A) edge (C)
		(D) edge[dashed] (A)
		(C) edge[dashed] (E)
			;
		\end{scope}
		
	\end{scope}
	
	\begin{scope}[shift=({8,1})]
			\begin{scope}
			\node at (-1,-0.3) {$R_2$};
			\node (R)[parallel] {$p$}
			child[sibling distance=15mm] {node [series] {$s$}
			child[sibling distance=7mm] {node [leaf] {$e$}}
			child[sibling distance=7mm] {node [leaf] {$e$}}}
			child[sibling distance=15mm] {node [series] {$s$}
			child[sibling distance=7mm] {node [leaf] {$e$}}
			child[sibling distance=7mm] {node [leaf] {$e$}}}
			;\draw[dashed] (R) -- (1.3,-.8);
			\end{scope}
		
		\begin{scope}[shift=({0,-3})]
			\path   
			node (A) [vertex=above $s_1$]     at ( -1, 0) {} 
			node (B) [vertex=above $v_1$]     at ( 0, .5) {} 
			node (C) [vertex=above $t_1$]     at ( 1, 0) {} 
			node (D) [vertex=below $v_2$]     at ( 0, -.5) {} 
				node (F)     at ( -0.7, -0.5) {} 
			node (E)      at ( 0.7,-0.5) {} 
			
			(A) edge (B)
			(B) edge (C)
			(A) edge (D)
			(D) edge (C)
				(F) edge[dashed] (A)
			(C) edge[dashed] (E)
			;
		\end{scope}
		
	\end{scope}

	\begin{scope}[shift=({7,-3})]
	\begin{scope}
		\node at (-1,-0.3) {$R_{2,1}$};
		\node (R)[series] {$s$}
		child[sibling distance=15mm] {node [parallel] {$p$}
		child[sibling distance=15mm] {node [series] {$s$}
			child[sibling distance=7mm] {node [leaf] {$e$}}
			child[sibling distance=7mm] {node [leaf] {$e$}}}
		child[sibling distance=15mm] {node [series] {$s$}
			child[sibling distance=7mm] {node [leaf] {$e$}}
			child[sibling distance=7mm] {node [leaf] {$e$}}}}
		child {node [leaf] {$e$}}
		;\draw[dashed] (R) -- (1.3,-.8);
		\draw[dashed] (R) -- (-1.3,-.8);
	\end{scope}
	
	\begin{scope}[shift=({-.5,-4})]
		\path   
		node (A) [vertex=above $s_1$]     at ( -1, 0) {} 
		node (B) [vertex=above $v_1$]     at ( 0, .5) {} 
		node (C) [vertex=above $t_1$]     at ( 1, 0) {} 
		node (D) [vertex=below $v_2$]     at ( 0, -.5) {} 
		node (E) [vertex=above $t_2$]     at ( 2, 0) {} 
			node (F)     at ( -1.7, -0) {} 
		node (G)      at ( 2.7, -0) {} 
		
		(A) edge (B)
		(B) edge (C)
		(A) edge (D)
		(D) edge (C)
		(C) edge (E)
			(F) edge[dashed] (A)
		(E) edge[dashed] (G)
		;
	\end{scope}
	
\end{scope}

\begin{scope}[shift=({1,-2.5})]
	\begin{scope}[shift=({0,-.5})]
		\node at (-1,-0.3) {$R_3$};
		\node (R)[parallel,sibling distance=15mm] {$p$}
		child[sibling distance=15mm] {node [series] {$s$}
			child[sibling distance=7mm] {node [leaf] {$e$}}
			child[sibling distance=7mm] {node [leaf] {$e$}}}
		child[sibling distance=15mm] {node [series] {$s$}
			child[sibling distance=7mm] {node [leaf] {$e$}}
			child[sibling distance=7mm] {node [leaf] {$e$}}}
			child[sibling distance=15mm] {node [series] {$s$}
				child[sibling distance=7mm] {node [leaf] {$e$}}
				child[sibling distance=7mm] {node [leaf] {$e$}}}
		;\draw[dashed] (R) -- (2,-.8);
	\end{scope}
		\begin{scope}[shift=({0,-4})]
		\path   
		node (A) [vertex=above $s_1$]     at ( -1, 0) {} 
		node (B) [vertex=above $v_1$]     at ( 0, .8) {} 
		node (C) [vertex=above $t_1$]     at ( 1, 0) {} 
		node (D) [vertex=below $v_3$]     at ( 0, -.8) {} 
		node (E) [vertex=above $v_2$]     at ( 0, 0) {} 
		node (F)     at ( -0.5, -1) {} 
		node (G)     at ( 0.5, -1) {} 
		
		(A) edge (B)
		(B) edge (C)
		(A) edge (D)
		(D) edge (C)
		(C) edge (E)
		(A) edge (E)
		(A) edge[dashed] (F)
		(C) edge[dashed] (G)
		;
	\end{scope}
	\end{scope}	
	
		\begin{scope}[shift=({4.5,-8})]
		\begin{scope}
			\node at (-1,-0.1) {$R_{3,2}$};
		\node (R)[series] {$s$}
		child[sibling distance=40mm] {node [parallel] {$p$}
			child[sibling distance=15mm] {node [series] {$s$}
				child[sibling distance=7mm] {node [leaf] {$e$}}
				child[sibling distance=7mm] {node [leaf] {$e$}}}
			child[sibling distance=15mm] {node [series] {$s$}
				child[sibling distance=7mm] {node [leaf] {$e$}}
				child[sibling distance=7mm] {node [leaf] {$e$}}}
			child[sibling distance=15mm] {node [series] {$s$}
				child[sibling distance=7mm] {node [leaf] {$e$}}
				child[sibling distance=7mm] {node [leaf] {$e$}}}
			}
		child[sibling distance=35mm] {node [parallel] {$p$}
		child[sibling distance=15mm] {node [series] {$s$}
			child[sibling distance=7mm] {node [leaf] {$e$}}
			child[sibling distance=7mm] {node [leaf] {$e$}}}
		child[sibling distance=15mm] {node [series] {$s$}
			child[sibling distance=7mm] {node [leaf] {$e$}}
			child[sibling distance=7mm] {node [leaf] {$e$}}}
	}
		;\draw[dashed] (R) -- (3.5,-.8);
		\draw[dashed] (R) -- (-3.5,-.8);
	\end{scope}

\begin{scope}[shift=({-1,-4})]
			\path   
			node (A) [vertex=above $s_1$]     at ( -1, 0) {} 
			node (B) [vertex=above $v_1$]     at ( 0, .8) {} 
			node (C) [vertex=above $t_1$,vertex=below $s'_1$]     at ( 1, 0) {} 
			node (D) [vertex=below $v_3$]     at ( 0, -.8) {} 
			node (E) [vertex=above $v_2$]     at ( 0, 0) {} 
			node (H) [vertex=above $v'_1$]     at ( 2, 0.5) {} 
			node (I) [vertex=below $v'_2$]     at ( 2, -0.5) {} 
			node (J) [vertex=below $t'_1$]     at ( 3, 0) {} 
				
			node (F)     at ( -1.8, 0) {} 
			node (G)     at ( 3.8, 0) {} 
				
			(A) edge (B)
			(B) edge (C)
			(A) edge (D)
			(D) edge (C)
			(C) edge (E)
			(A) edge (E)
			(C) edge (H)
			(C) edge (I)
			(J) edge (H)
			(J) edge (I)
			(A) edge[dashed] (F)
			(J) edge[dashed] (G)
			;		
		\end{scope}	
		
	\end{scope}	
\end{tikzpicture}

			\caption{Substructures $P_2, Q_2, R_2, R_3, R_{2,1}, R_{3,2} $ and the corresponding subgraphs.}\label{fig:substr}
		\end{center}
	\end{figure}
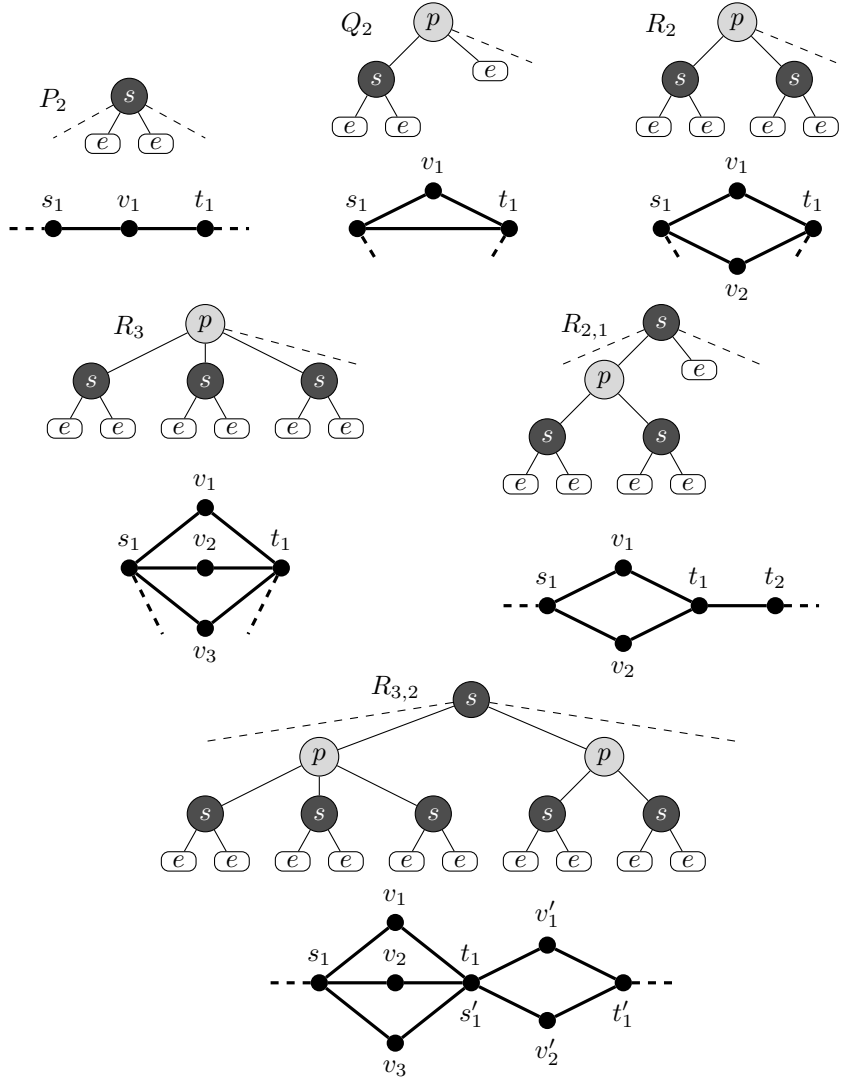

	\begin{defn} Let $P_2$ denote the substructure consisting of a series connection of two edges, and let $Q_2$ denote the parallel join of $P_2$ with a single edge.
	For $i \geq 2$, define $R_i$ to be the substructure given by the parallel join of $i$ copies of $P_2$.
	Then, for $i, j \geq 1$, let $R_{i,j}$ denote the series join of $R_i$ and $R_j$, where we define $R_1$ to be a single edge.
	(In particular, $R_{1,1} = P_2$.)
	See Fig.~\ref{fig:substr} for illustrations of these substructures. 
		\end{defn}
	
	{In Fig.~\ref{sp-tree} the reader can see examples of these substructures formed by the following edges. $P_2=e_{12},e_{13}$; $Q_2=e_{27},e_{28}, e_{29}$; $R_2=e_{18},e_{19}, e_{20}, e_{21}$;  $R_{2,1}=e_{17},\ldots, e_{21}$.  
	The graph contains no \(R_3\) substructure, but the edges \(e_2,\ldots,e_8\) form one after replacing $e_5,e_6$ with a single edge. 	Also, there is no $R_{3,2}$ in this example.}


		\begin{lem}\label{lem:middle}
		Let $G$ be an sp-graph containing the substructure $R_2$
		or $R_3$, and suppose that $\tau(G)\le\frac12$.
		If $S$ is a tough set in $G$, then the middle vertices
		of these substructures do not belong to $S$.
	\end{lem}

	\begin{proof}
		We present the proof for $R_2$, the argument is almost the same for $R_3$. 
		It is easy to see that if at least one of $s_1$ and $t_1$ belongs to $S$ as well as at least one of $v_1$ and $v_2$, then by setting $S'=S-\{v_1,v_2\}$ we have $|S'|<|S|$ and $c(G-S')\ge c(G-S)\ge 2$. Thus $S$ cannot be a tough set. Hence, we can assume that $s_1,t_1\notin S$.
		
		Now suppose that $v_1\in S$ and $v_2\notin S$. In this case, $s_1$ and $t_1$ belong to the same component of $G-S$ since the path $s_1v_2t_1$ connects them. By setting $S'=S-v_1$ we have $|S'|<|S|$ and $c(G-S')=c(G-S)$. Thus $S$ cannot be a tough set, again.
		
		Finally, suppose $v_1,v_2\in S$ and let $S'=S-\{v_1,v_2\}\cup t_1$. Clearly,  $|S'|<|S|$. Now consider the components of $G-S'$; how are these different from the components of $G-S$? We can see that $v_1$ and $v_2$ will belong to the same component as $s_1$. Since any path of $G$ between different components of $G-S$ containing $v_1$ or $v_2$ also contains $t_1$, it is impossible that two vertices in different components of $G-S$ are connected with a path in $G-S'$. If $t_1$ is an isolated vertex in $G-S$ then the number of components may decrease by 1. On the other hand, adding $t_1$ to the cutset may separate some other components. Therefore, $c(G-S')\ge c(G-S)-1$. Since $\tau(G)\le\frac{1}{2}$, we have $c(G-S')\ge c(G-S)-1\ge 2|S|-1\ge 3$, which implies that $S'$ is a cutset giving a better ratio than $S$, so $S$ cannot be a tough set.
	\end{proof}
	
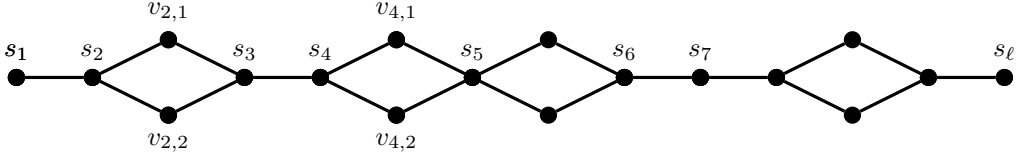
\begin{figure}[h]
		\begin{center}
%
%
%

	\begin{tikzpicture}[scale=1,
	every edge/.style = {draw,very thick},
	vertex/.style args = {#1 #2}{circle, 
		draw, 
		thick,
		minimum size=2pt,inner sep=-2pt,
		fill=black,
		label=#1:#2},
			leaf/.style={rectangle, draw=black, rounded corners=1mm, 
			text centered, anchor=north, text=black,text width=1em,inner sep=1.5pt},
		parallel/.style={circle, draw=black, fill=gray!30,
			text centered, anchor=north, text=black,text width=1em,inner sep=1pt},
		non/.style={circle, 
				text centered, anchor=north, text=black,text width=1em,inner sep=1pt},
		series/.style={circle, draw=black, fill=black!70, 
			text centered, anchor=north, text=white, text width=1em,inner sep=1pt},
		level distance=0.25cm, growth parent anchor=south,level/.style={sibling distance=10mm}]
		
	
	\begin{scope}[shift=({0,0})]
		\foreach \x in {1,4,6,10}{
		\begin{scope}[shift=({\x,0})]
			\path   
			node (A)  [vertex=above $$]    at ( 0, 0) {} 
			node (B)  [vertex=above $$]   at ( 1, .5) {} 
			node (C) [vertex=above $$]     at ( 2, 0) {} 
			node (D) [vertex=below $$]     at ( 1, -.5) {}

			(A) edge (B)
			(B) edge (C)
			(A) edge (D)
			(D) edge (C)
			
			;
		\end{scope}
	}
		\foreach \x in {0,3,8,9,12}{
		\begin{scope}[shift=({\x,0})]
			\path   
			node (A)  [vertex=above $$]    at ( 0, 0) {} 
			node (B)  [vertex=above $$]   at ( 1, 0) {}

			(A) edge (B)

			;
		\end{scope}
	}
	
	\path 
	node  [vertex=above $s_1$] at (0,0) {}
	node  [vertex=above $s_2$] at (1,0) {}
	node  [vertex=above $s_3$] at (3,0) {}
	node  [vertex=above $s_4$] at (4,0) {}
	node  [vertex=above $s_5$] at (6,0) {}
	node  [vertex=above $s_6$] at (8,0) {}
	node  [vertex=above $s_7$] at (9,0) {}
	node  [vertex=above $s_1$] at (0,0) {}
	node  [vertex=above $s_{\ell}$] at (13,0) {}
	node  [vertex=above $v_{2,1}$] at (2,.5) {}
	node  [vertex=below $v_{2,2}$] at (2,-0.5) {}
	node  [vertex=above $v_{4,1}$] at (5,.5) {}
	node  [vertex=below $v_{4,2}$] at (5,-0.5) {}
	;

	\end{scope}

\end{tikzpicture}

			\caption{A necklace graph.}\label{fig:necklace}
		\end{center}
	\end{figure}

	\begin{defn}A \emph{necklace graph} is a series join of edges and $R_2$
		substructures, arranged so that both the first and last joins are single edges. (See Fig.~\ref{fig:necklace}.) 
	\end{defn}
	
	A path is a special case of necklace graphs.
	
	\begin{lem}\label{necklace}
		The necklace graphs are minimally $\frac{1}{2}$-tough.
	\end{lem}
	\begin{proof}
		The toughness is $\leq \frac{1}{2}$ since there is a cut-vertex, for example, the joining vertex at the end of the first edge, namely \(s_2\).
		Suppose that $\tau (G)< \frac{1}{2}$ and let $S$ be a tough set.
		
		By Lemma~\ref{lem:middle},  $S$ does not contain either of the vertices \(v_1\) and \(v_2\) of the $R_2$ substructures. Therefore, $S$ consists of some of the joining vertices $s_1,\ldots,s_{\ell}$ (see Fig.~\ref{fig:necklace} for the notation).
		
		Suppose $s_i\in S$, but $s_j\notin S$ for all $1\le j<i$. Then in $G-S$ there will be one component to the left of $s_i$, this contains $s_1,\ldots,s_{i-1}$. A symmetric situation occurs on the right.
		
		If $s_i,s_{i+1}\in S$, then there are two components between  $s_i,s_{i+1}$, the isolated vertices $v_{i,1}$ and $v_{i,2}$. 	If $s_i,s_{j}\in S$ with $j>i+1$ and $s_{i+1},\ldots, s_{j-1}\notin S$, then  a single component $G-S$  contains all its vertices between  $s_i$ and $s_{j}$. So to summarize, $G-S$ contains at most two components between two consecutive $s_i$ vertices in $S$, and there may be one additional component on each end.
		 Hence, we obtain
		$$\frac{|S|}{c(G-S)}\geq \frac{k}{(2k-2)+2}=\frac{1}{2}.$$
				
		Now we show that for any edge $\tau(G-e)< \tau(G).$
		If the edge is a cut-edge, then $\tau (G-e)=0 < \tau(G)= \frac{1}{2}.$

		Otherwise, without loss of generality, we can assume that $e=(s_i,v_{i,1})$ for some $R_2$ subgraph. In this case, $c((G-e)-s_{i+1})=3$, so the toughness is at most $\frac{1}{3}< \frac{1}{2}$.
	\end{proof}

	The next lemma is crucial in the following arguments. We will only use it for $k=2$, but it might be interesting in this more general form. Note that the claim is not necessarily true if the toughness is not in the form $\frac{1}{k}$.
	
	\begin{lem} \label{k}
		Let $G$ be a graph with $\tau(G)=\frac{1}{k}$, for some integer $k \geq2$. Let $e$ be an edge so that $\tau(G-e)< \tau(G)$ and $S$ be a tough set in $G-e$. Now $S$ is a tough set in $G$ as well.
	\end{lem}

	\begin{proof}
		If $\tau(G)=\frac{1}{k}$ then $\tau(G-e)=\frac{|S|}{c((G-e)-S)} < \frac{1}{k}$ by the choice of $e$. This implies that $ c((G-e)-S)>2.$
	
		Adding $e$ back to $(G-e)-S$ either leaves the number of components unchanged or decreases it by 1.
		 Therefore, $1< c((G-e)-S)-1 \leq c(G-S) \leq c((G-e)-S)$,
			which shows that $S$ is a cutset in $G$.
		
			Now observe that
		$$\tau(G)= \frac{1}{k} \leq \frac{|S|}{c(G-S)} \leq \frac{|S|}{c((G-e)-S)-1}$$ 
		while from the assumption on $\tau(G - e)$, we also have
		$$ \frac{|S|}{c((G-e)-S)}< \frac{1}{k} \leq \frac{|S|}{c((G-e)-S)-1}.$$ 
		These inequalities imply that
		$ c((G-e)-S)-1 \leq k |S| < c((G-e)-S)$, 	so $k|S|$ lies between two consecutive integers and is therefore itself an integer. Hence, we conclude that 	
		$c((G - e) - S) - 1 = k|S|$.
		
		Putting everything together, we have
		$$\frac{1}{k} \leq \frac{|S|}{c(G-S)} \leq \frac{|S|}{c((G-e)-S)-1}= \frac{1}{k},$$ 	which shows that $\frac{|S|}{c(G - S)} = \frac{1}{k}$. Therefore, $S$ is a tough set in $G$.
		\end{proof}

	Next, we show that minimally $\frac{1}{2}$-tough graphs cannot contain certain substructures. 

  \begin{lem} \label{lem:R4}
  	If  an sp-graph $G$ contains $R_4$ as a substructure, then it is not minimally $\frac{1}{2}$-tough.
  \end{lem}
  
  \begin{proof}
  	First, suppose that $G=R_4$ or $G$ is the parallel join of $R_4$ and an edge. It is easy to verify that in both cases, $\tau(G)=\frac{1}{2}=\tau(G-e)$ holds for any edge $e$ of the graph. Hence $G$ is not minimally $\frac{1}{2}$-tough.
  	
  	In any other case, $G$ contains at least one vertex apart from the vertices of $R_4$. Setting $S=\{s_1,t_1\}$ gives $c(G-S)\ge 4+1=5$, so $\tau(G)\le\frac{2}{5}$, therefore $G$ is not even $\frac{1}{2}$-tough.
  \end{proof}
  	
  	This implies that we can assume that $G$ does not contain the substructure $R_4$. Note that this also implies that $G$ does not contain any $R_i$ with $i\ge 4$.
  	
  \begin{lem}\label{lem:Q2}
  		If  a minimally $\frac{1}{2}$-tough sp-graph $G$ contains the substructure $Q_2$, then $G$ is not minimally $\frac{1}{2}$-tough.
  \end{lem}
  \begin{proof}
  		Let $e=s_1t_1$, we show that $\tau(G-e)=\frac{1}{2}$. Suppose on the contrary that $\tau(G-e)<\frac{1}{2}$, then there exists a minimum-size tough set $S$ for $G-e$. By Claim~\ref{claim} we know that  $s_1,t_1 \notin S$. By Lemma~\ref{l3}, either $S= \{v_1\}$ or $S$ is disjoint from $Q_2.$
  		The first is impossible since  Lemma~\ref{k} implies that $S$ must be a tough set in $G$ as well, but it is not even a cutset. In the second case, $c((G-e)-S)= c(G-S)$ since $s_1$ and $t_1$ is connected by the path $s_1 v_1 t_1$. Hence, $\tau(G-e)<\frac{1}{2}$ cannot hold, a contradiction.
  \end{proof} 
  
  \begin{lem} \label{lem:R33}
  	If a minimally $\frac{1}{2}$-tough sp-graph $G$ contains $R_{3,3}$ as a substructure, then it is not minimally $\frac{1}{2}$-tough.
  \end{lem}
  
  \begin{proof}
  	First, suppose that $G=R_{3,3}$ or $G$ is the parallel join of $R_{3,3}$ and an edge. It is easy to verify that in both cases, $\tau(G)=\frac{1}{2}=\tau(G-e)$ holds for $e=v_1t_1$. Hence $G$ is not minimally $\frac{1}{2}$-tough. (Notice that
  	for $e=s_1v_1$ the vertex set $\{s'_1,t'_1\}$ is a tough set in $G-e$, showing that $\tau(G-e)=\frac{2}{5}$.)
  	
  	In any other case, $G$ contains at least one vertex apart from the vertices of $R_{3,3}$. Setting $S=\{s_1,t_1,t'_1\}$ gives $c(G-S)\ge 6+1=7$, so $\tau(G)\le\frac{3}{7}$, therefore $G$ is not even $\frac{1}{2}$-tough.
  \end{proof}
  
  \begin{lem} \label{lem:R32}
  	If   a minimally $\frac{1}{2}$-tough sp-graph $G$ contains $R_{3,2}$ as a substructure, then it is not minimally $\frac{1}{2}$-tough.
  \end{lem}
  
  \begin{proof}
  	First, suppose that $G=R_{3,2}$ or $G$ is the parallel join of $R_{3,2}$ and an edge. It is easy to verify that in both cases, $\tau(G)=\frac{1}{2}=\tau(G-e)$ holds for $e=v_1t_1$. Hence $G$ is not minimally $\frac{1}{2}$-tough. 
  	
  	In any other case, $G$ contains at least one vertex apart from the vertices of $R_{3,2}$. 
  	Let $e=v_1t_1$, we show that $\tau(G-e)=\tau(G)=\frac{1}{2}$, which proves our claim.
  	
  	Suppose to the contrary that  $\tau(G-e)<\tau(G)$ and let $S$ be a tough set of $G-e$. Lemma~\ref{k} implies that $S$ is a tough set of $G$ as well.  By Lemma~\ref{lem:middle}, none of the vertices $v_i$ or $v'_i$ belong to $S$. By Claim~\ref{claim} we have  $v_1\notin S$ and $t_1\notin S$. This implies that $s_1\in S$, otherwise $v_1$ and $t_1$ belong to the same component of $G-e$, which contradicts the assumption that $S$ is a tough set.
  	
  	Let $S'=S\cup\{t_1\}$. It is clear that in $G-S$, the vertices $v_1,v_2,v_3, t_1,v'_1$, and $v'_2$ all belong to the same component. Now, if $t_1$ is deleted from $G-S$, then $v_1,v_2,v_3$ become isolated vertices, and there will be a component that contains $v'_1$. This component will contain $v'_2$ as well if $t'_1\notin S$; otherwise, $v'_1$ and $v'_2$ also become isolated vertices.  In either way, we obtain
  	$c(G-S')\ge c(G-S)+3$. By the Mediant Inequality (Lemma~\ref{mediant}) we have 
  	$$\frac{1}{3}< \frac{|S|+1}{c(G-S)+3}< \frac{|S|}{c(G-S)}= \frac{1}{2}.$$
  	Thus 
  	$$\tau(G) \leq \frac{|S'|}{c(G-S')}\leq \frac{|S|+1}{c(G-S)+3}<  \frac{1}{2},$$ which contradicts that $\tau (G) = \frac{1}{2}$, since $S'$ is clearly a cutset.
  \end{proof}

  \begin{lem}\label{lem:R31}
  	If   a minimally $\frac{1}{2}$-tough sp-graph $G$ contains $R_{3,1}$ as a substructure, then it is not minimally $\frac{1}{2}$-tough.
  \end{lem}
  
\begin{proof}
	First, suppose that $G=R_{3,1}$ or $G$ is the parallel join of $R_{3,1}$ and an edge. It is easy to verify that in both cases, $\tau(G)=\frac{1}{2}=\tau(G-e)$ holds for $e=v_1t_1$. Hence $G$ is not minimally $\frac{1}{2}$-tough. 
	
	Now suppose that the parent of  $R_{3,1}$ is the root of $G$, so  $G$ is a series join of $R_{3,1}$ and some other subgraphs. If there is a series join on both sides (i.e., $s_1$ and $t_2$ are both join vertices) then $S=\{s_1,t_1\}$ is clearly a cutset with $c(G-S)=5$, a contradiction.
	If there is no series join at $s_1$, then let $e=v_1t_1$ and $S$ be a tough set of $G-e$. By Lemma~\ref{k} we know that $S$ is a tough set in $G$ as well. By Claim~\ref{claim} we have $t_1,v_1 \notin S$. Lemma~\ref{lem:middle} implies that  $v_2,v_3\notin S$. The path $t_1 v_2 s_1 v_1$ connects $t_1$ and $v_1$ in $G-e$, therefore $s_1 \in S$ must hold, otherwise $t_1$ and $v_1$ belong to the same component of $(G-e)-S$, a contradiction. It is clear that in this case $S\not=\{s_1\}$ since it would not be a cutset in $G$. Thus, it is easy to verify that $S'=S-s_1$ is still a cutset in $G$, and that $c(G-S')=c(G-S)$ holds. This contradicts the assumption that $S$ is a tough set in $G$.
	
	The above argument shows that, if the root of \(R_{3,1}\) has a parent, that parent must be a parallel node. Let $e=s_1v_1$, we show that $\tau (G-e) \geq \tau (G)$, which implies that it is not minimally tough. 
	
	Let $S$ be a tough set in $G-e$. By Claim \ref{claim} we have $s_1,v_1 \notin S$. Lemmas~\ref{k} and \ref{lem:middle} imply that  $v_2,v_3\notin S$.  The path $s_1 v_2 t_1 v_1$ connects $s_1$ and $v_1$ in $G-e$, therefore $t_1 \in S$ must hold, otherwise $s_1$ and $v_1$ belong to the same component of $(G-e)-S$, a contradiction. 
	
	Thus, the remaining case is when $t_1 \in S$ and $s_1,v_1,v_2,v_3 \notin S$.	
	Let $S'=S-t_1$, hence $|S'|=|S|-1$. We claim that $c(G-S') \geq c(G-S)-1$.
	Consider the components of $G-S$. Three of the neighbors of $t_1$, namely $v_1,v_2,v_3$ are  in the same component, since $s_1,v_1,v_2,v_3 \notin S$. $t_2$ may be in a different component, but restoring $t_1$  may only merge these two components, so the number of components may only decrease by $1$.
	
	We also claim that $c(G-S')>1$, i.e., $ S'$ is a cutset. Lemma~\ref{k}
	implies that $S$ is a tough set also in $G$. Since the parent of the root of $ R_{3,1}$ is a parallel node,  there is a path between $s_1$ and $t_2$ that avoids the substructure. This implies that $G-t_1$ is connected, so 	$S \not= \{t_1\}$, and hence $|S|\geq 2$. Since $\tau(G)= \frac{1}{2}$, this gives that $c(G-S) \geq 4$ and $c(G-S') \geq 3 >1$.
	
	Using the above claims, we obtain $$\tau(G)\leq \frac{|S'|}{c(G-S')} \leq \frac{|S|-1}{c(G-S)-1}.$$
	
	However, $\frac{1}{2}=\frac{|S|}{c(G-S)}$  is the mediant of $\frac{|S|-1}{c(G-S)-1}$ and $\frac{1}{1}$. Thus, the Mediant Inequality implies that
	$$\tau(G) \leq \frac{|S'|}{c(G-S')}\leq \frac{|S|-1}{c(G-S)-1}< \frac{|S|}{c(G-S)}= \frac{1}{2}<1,$$ which contradicts the assumption that $\tau (G) = \frac{1}{2}$.
\end{proof}

  \begin{lem}\label{lem:R21}
  	If a minimally $\frac{1}{2}$-tough sp-graph $G$ contains the substructure $R_{2,1}$ then the root of $R_{2,1}$ is the root of $G$ and $s_1$ is a cut-vertex.
  \end{lem}
  
  Notice that a minimally $\frac{1}{2}$-tough sp-graph may contain  $R_{2,1}$ substructures, some of the necklace graphs have this property.
  
  \begin{proof}
  The proof is essentially the same as the proof of Lemma~\ref{lem:R31}. The main difference is that if $s_1$ is a cut-vertex, then the above argument does not work, hence we have the exceptional case.
  \end{proof}

  So far, we have shown that a minimally $\frac{1}{2}$-tough sp-graph cannot contain certain substructures. This also holds for the reduced $\frac{1}{2}$-tough sp-graphs, which implies that their structure is fairly simple. Now we turn our attention to the height of the sp-tree of the reduced graph. Our goal is to prove that the height is at most 3. In fact, we show that it cannot have a subtree of height 4 or more. To achieve this, we characterize the possible subtrees of height $0$--$3$. (Note that now we are considering subtrees, not substructures.)
  
  \begin{lem}\label{heights}
  	If $G$ is a minimally $\frac{1}{2}$-tough reduced sp-graph, then the height of its sp-tree can only be 1 or 3. If the height is $1$, then $G=P_3$, otherwise, $G$ is a necklace graph different from $P_3$.
  \end{lem}
  
  \begin{proof}
  Let us recall a few facts about sp-trees:
  \begin{itemize}
  	\item  The height of a tree is the number of edges on the longest path from the root to a leaf.
  	\item Every non-leaf node has at least two children.
  	\item If the height of a tree is $k$, then at least one child of the root is the root of a subtree of height $k-1.$
  	\item The parent of a series node is a parallel node and vice versa.
  \end{itemize}
  
  If $v$ is a node of the sp-tree with children $w_1,\ldots,w_k$, then we call the subtrees with roots $w_1,\ldots,w_k$ the \emph{direct subtrees} of $v$.
  
  \subsubsection*{Height $0$:} The sp-tree has only one node, the root. This node is a leaf representing an edge of the graph.
  Since the toughness of a single edge is $\infty$, this cannot be the whole graph, it can only be a proper subtree.
  
  \subsubsection*{Height $1$:} The root cannot be a parallel node, since this would give parallel edges, contradicting Lemma~\ref{l5}. Thus, the root is a series node. It must have at least two children, all of whose children are leaves, each representing an edge. On the other hand, it cannot have more than two children that are edges, because the graph is reduced. This subtree corresponds to $P_3$ (i.e., the path on 3 vertices with 2 edges). The graph $P_3$ itself is clearly minimally $\frac{1}{2}$-tough.
  
  \subsubsection*{Height $2$:} At least one direct subtree must have  height $1$, such a subtree can only be a  $P_3$.
  
  Since the root of $P_3$ is a series node, the root of the height $2$ tree must be a parallel node. This node cannot have more than one edge as a child, as parallel edges are not allowed.  Moreover, if it has both an edge child and a $P_3$ child, then a $Q_2$ substructure would arise, contradicting Lemma~\ref{lem:Q2}.
  
  Thus, all the children of the root must be $P_3$ graphs. Having four or more such children would yield an $R_4$ substructure, which is excluded by Lemma~\ref{lem:R4}.	Hence, the subtree must have exactly two or three children, each being $P_3$, so the subtree is either $R_2$ or $R_3$.
  
   Since $R_2$ has toughness 1 and $R_3$ has toughness $\frac{2}{3}$, neither can be the entire graph, they can only occur as proper subtrees.
  
  \subsubsection*{Height $3$:} At least one direct subtree must have height $2$, and such a subtree can only be an $R_2$ or $R_3$.  Since their roots are parallel nodes, the root of the height $3$ tree must be a series node. This implies that no direct subtree can have a series node as its root, which rules out subtrees of height 1.  Thus, all children of the root are either $R_2$ or $R_3$ subtrees or single edges. That is, the entire graph is a series join of $R_2$, $R_3$ substructures, and edges, in some order.
  
   We now show that if this subtree is the whole graph, then it is minimally $\frac{1}{2}$-tough if and only if it contains no $R_3$ substructure, and both the first and last graphs in the series composition must be edges, not $R_2$ substructures.
  
   Let $S = \{s_1, t_1\}$, where $s_1$ and $t_1$ are the terminals connecting the $R_3$ substructure to the rest of the graph, so $R_3$ is not the first or last substructure. Then $S$ forms a cutset of size 2, and the removal of $S$ results in 5 components: $c(G - S) = 5$. Therefore, $\tau(G) \le \frac{2}{5}$, a contradiction. 
  
   Next, suppose that the first direct subtree is an $R_2$ or $R_3$ with terminals $s_1$ and $t_1$, where $s_1$ is the initial terminal of the graph (i.e., not a cut-vertex) and $t_1$ is an internal terminal (i.e., a cut-vertex). Using arguments similar to those above, one can show that removing the edge $e = t_1v_1$ does not decrease the toughness, implying that $G$ is not minimally tough.
  
  Hence, the graph must be a necklace graph.

	\subsubsection*{Height $4$ or more:}  If the tree has height at least 4, then it must contain at least one subtree of height exactly 4. If there is more than one such subtree, we select the one that is at maximum distance from the root of the entire tree. Let $T_4$ denote this subtree, and let $r_4$ be its root. This choice guarantees that if $r_4$ has siblings, then the subtrees rooted at those siblings have height at most 4.
	
	One of the direct subtrees of $r_4$ must have height $3$, and hence is a series join of $R_2$, $R_3$ substructures, and edges. Since its root is a series node, $r_4$ must be a parallel node. Consequently, $r_4$ cannot have a direct subtree of height $2$.
	
	Now suppose that the height $3$ subtree contains an $R_3$ substructure. Consider its neighboring subtree in the series join. This neighbor must be either an edge, an $R_2$, or another $R_3$, resulting in an $R_{3,1}$, $R_{3,2}$, or $R_{3,3}$ substructure, respectively. However, Lemmas~\ref{lem:R31}, \ref{lem:R32}, and \ref{lem:R33} rule out these configurations, leading to contradictions in all cases. Therefore, the height $3$ subtree cannot contain an $R_3$ substructure.

	Now suppose that the height $3$ subtree contains at least one $R_2$ and at least one edge among the children of the series node. Then, necessarily, there exists an $R_2$ substructure and an edge that are adjacent in the series join, forming an $R_{2,1}$ substructure. In this configuration, the common ancestor of these substructures is the parallel node $r_4$. Therefore, by Lemma~\ref{lem:R21}, such a configuration is not allowed, a contradiction.
	
	This implies that any height $3$ subtree must be a series join of two or more $R_2$ graphs. The other direct subtrees of $r_4$ may have height 1 or 0; that is, they may be $P_3$ subgraphs or single edges.
	
	Let us now analyze the possible configurations. In each case, the argument differs depending on whether $T_4$ is a proper subtree or whether it is the whole graph. Subcases (a) cover the first case, while subcases (b) cover the second.
	
	We define \emph{a bracelet of length $\ell$}, denoted by $B_\ell$, as a series join of $\ell$ copies of the $R_2$ graph. Let $s_1, s_2, \ldots, s_{\ell+1}$ denote the joining vertices of the consecutive $R_2$ subgraphs, in natural order. The terminals of the bracelet are $s_1$ and $s_{\ell+1}$.
	Observe that removing the set $\{s_1, s_2, \ldots, s_{\ell+1}\}$ from the bracelet leaves $2\ell$ components, each of which is an isolated vertex.

  \subsubsection*{Case 1(a): One child of $r_4$ is $B_{\ell}$, another child is an edge $e$, and $r_4$ is not the root.} 
 
  Since $e$ is a jump-edge and its parent is not the root, Lemma~\ref{j1} implies that the graph is not minimally $\frac{1}{2}$-tough.

  \subsubsection*{Case 1(b): One child of $r_4$ is $B_{\ell}$, another child is an edge $e$, and $r_4$ is the root.}
  
  If $r_4$ has at least three children, then $e$ is again a jump-edge. By Lemma~\ref{j1}, this implies that the graph is not minimally $\frac{1}{2}$-tough.
  
  Now suppose that $r_4$ has exactly two children. Then the graph is a parallel join of a bracelet $B_\ell$ and a single edge $e$. We claim that in this case $\tau(G) > \frac{1}{2}$.
  
  Let $S$ be a minimum-size tough set of $G$. By Lemma~\ref{lem:middle}, $S \subseteq \{s_1, s_2, \ldots, s_{\ell+1}\}$. If $S = \{s_1, s_2, \ldots, s_{\ell+1}\}$, then
  $\tau(G) = \frac{\ell+1}{2\ell} > \frac{1}{2}$.
   Suppose instead that two vertices $s_p, s_q \notin S$ are such that $G - \{s_p, s_q\}$ has two nonempty components, each intersecting $S$. This contradicts Lemma~\ref{l3}, which rules out such a configuration in a minimum-size tough set. Therefore, the only possible forms of $S$ are either:
   $S = \{s_i, s_{i+1}, \ldots, s_j\}$, for some $1 \le i < j \le \ell+1$; or
   $S = \{ s_1, \ldots, s_i, s_j, s_{j+1}, \ldots, s_{\ell+1}\}$, for some  indices $1 \le i < j \le \ell+1$.
   In the first case,
  $$
  \tau(G) = \frac{j - i + 1}{2(j - i)+1} > \frac{1}{2},
  $$
  and in the second case,
  $$
  \tau(G) = \frac{(\ell + 1 - j + 1 + i)}{2(\ell + 1 - j) + 2(i - 1)+1} > \frac{1}{2}.
  $$
   In both cases, we conclude that $\tau(G) > \frac{1}{2}$, so the graph cannot be minimally $\frac{1}{2}$-tough.

  \subsubsection*{Case 2(a): One child of $r_4$ is $B_{\ell}$, another child is $B_{k}$, and $r_4$ is not the root.}
  
  Let $\{s_1, \ldots, s_{\ell+1}\}$ and $\{s'_1, \ldots, s'_{k+1}\}$ denote the joining vertices of the bracelets $B_\ell$ and $B_k$, respectively, chosen so that $s_1 = s'_1$ and $s_{\ell+1} = s'_{k+1}$. Define the set $S = \{s_1, \ldots, s_{\ell+1}\} \cup \{s'_1, \ldots, s'_{k+1}\}$.
  
  In the graph $G - S$, all the middle vertices of both bracelets become isolated, contributing $2\ell + 2k$ components. Moreover, since $r_4$ is not the root, there is at least one additional component in $G - S$. Therefore,
   $$
  \tau(G) \le \frac{|S|}{c(G - S)} \le \frac{\ell + k}{2\ell + 2k + 1} < \frac{1}{2},
  $$
    which contradicts the assumption that $\tau(G) \ge \frac{1}{2}$.
 
  \subsubsection*{Case 2(b): One child of $r_4$ is $B_{\ell}$, another child is $B_{k}$,  and $r_4$ is the root.}

  Suppose first that $r_4$ has only these two children.
  In this case, the graph consists of a cyclic structure formed by $\ell + k$ copies of the $R_2$ substructures. Let the joining vertices be $\{s_1, s_2, \ldots, s_{\ell+k}\}$, indexed cyclically. Consider the edge $e = s_1v_1$, where $v_1$ is one of the middle vertices of the first $R_2$, which has terminals $s_1$ and $s_2$. (By symmetry, any such edge behaves similarly.) We claim that $\tau(G - e) = \frac{1}{2}$, leading to a contradiction.
  
  Suppose  that $\tau(G - e) < \frac{1}{2}$, and let $S$ be a minimum-size tough set of $G - e$.  By Lemma~\ref{k}, $S$ is a tough set in $G$ as well. Using similar reasoning as in the previous cases, we can assume that $s_1 \notin S$ and $s_2 \in S$. By Lemmas~\ref{lem:middle} and~\ref{l3}, this implies that $S = \{s_2, s_3, \ldots, s_i\}$ for some $3 \le i < \ell + k$.
  
  Now observe that removing $i - 1$ vertices from $G$ disconnects the graph into $2(i - 2) + 1$ components. Therefore,
   $$
  \tau(G)= \frac{i - 1}{2(i - 2) + 1} > \frac{1}{2},
  $$
  which contradicts the assumption that $\tau(G) = \frac{1}{2}$. Hence, $\tau(G - e) = \frac{1}{2}$, and the graph is not minimally $\frac{1}{2}$-tough.
  
  If $r_4$ has any additional edge child $e$, then $G-e$ is the above investigated cyclic structure of $R_2$ substructures. Since $\tau(G-e)=\frac12=\tau(G)$, the graph is not minimally tough. 
  If $r_4$ has any other additional child besides $B_{\ell}$ and $B_k$, then the argument in Case 2(a) gives the same contradiction. 
  
  \subsubsection*{Case 3(a): One child of $r_4$ is $B_{\ell}$, two other children are $P_3$, and $r_4$ is not the root.}
 
  Let $S = \{s_1, s_2, \ldots, s_{\ell+1}\}$. In the graph $G - S$, all the middle vertices of the bracelet $B_\ell$ become isolated, contributing $2\ell$ components. Additionally, each of the two $P_3$ subgraphs contributes an isolated middle vertex, adding 2 more components. Since $r_4$ is not the root, there is at least one additional component in $G - S$. Therefore,
   $$
  \tau(G) \le \frac{|S|}{c(G - S)} \le \frac{\ell + 1}{2\ell + 3} < \frac{1}{2},
  $$
   which contradicts the assumption that $G$ is minimally $\frac{1}{2}$-tough.

  \subsubsection*{Case 3(b): One child of $r_4$ is $B_{\ell}$, two other children are $P_3$, and $r_4$ is the root.}
  
  The two $P_3$ subgraphs together form an $R_2$, so the graph is a cyclic structure consisting of $\ell + 1$ copies of $R_2$. If $r_4$ has no additional children, then the argument from Case 2(b) applies and shows that the graph is not minimally $\frac{1}{2}$-tough. $r_4$ cannot have an additional edge child, since that would give a $Q_2$ substructure, so by Lemma~\ref{lem:Q2} it is not minimally $\frac{1}{2}$-tough. On the other hand, if $r_4$ has any other additional children, then it is straightforward to verify that $\tau(G) < \frac{1}{2}$,  a contradiction.

  \subsubsection*{Case 4(a):  $r_4$ has only two children $B_{\ell}$ and $P_3$, and $r_4$ is not the root.}
  
  This is the most complicated case, and we need to consider several subcases. Let $\{s_1,\ldots, s_{\ell+1}\}$  denote the joining vertices of $B_{\ell}$ and let $s_1, w, s_{\ell+1}$ denote the vertices of the $P_3$ substructure.
  
   Since $r_4$ is not the root, its parent is a series node which must have at least two children, therefore $r_4$ must have at least one sibling. Without loss of generality, we may assume that it is joined to $T_4$ at $s_{\ell+1}$, let us denote it by $r'_4$, and let $T'_4$ be its subtree. By the choice of $r_4$, we know that $T'_4$ has height at most 4. Since $r'_4$ is either an edge or a parallel node, $T'_4$ can only be the following subtrees: an edge, $R_2$, $R_3$, or a height 4 subtree. However, in the previous cases we already ruled out most of the height 4 subtrees, the only remaining possibility is  a parallel join of $B_{k}$ and $P_3$. Now we consider each of these subcases.
  
  \medskip
  \textit{Subcase (i): $T'_4$ is an edge.}
  Let $s_{\ell+1}t'_1$ be the edge  $T'_4$.
  If there is no other vertex or edge in the graph, then let $e=ws_{\ell+1}$. Thus $G-e$ is a necklace graph, therefore its toughness is $\frac{1}{2}$, a contradiction.
  If $G$ also contains the leap edge $f=s_1t'_1$ then $\tau(G-f)\ge \tau(G-f-e)\ge \frac{1}{2}$, a contradiction again.

   \textit{Subsubcase (i:$\alpha$): The parent $r_5$ of $r_4$ is the root.} This means that $r_5$ is a series node, all its children are parallel nodes or edges, one of the children is $r_4$, the next child is $r'_4$. There may be other children on either side of them.
   
  First, suppose that   $r'_4$ has a sibling joined to $t'_1$; there may be further siblings joined to it. Let $H$ denote the subgraph induced by all these siblings, so it is a sequence of series-joined substructures. Thus $s_{\ell+1}t'_1$ is a bridge connecting $T_4$ and $H$. Let $e=ws_{\ell+1}$ again, we show that $\tau(G-e)=\tau(G)\le\frac{1}{2}$. Suppose, to the contrary, that  $\tau(G-e)<\tau(G)$ and let $S$ be a tough set of $G-e$. By Lemma~\ref{l2} $S$ is either contained in $T_4$ or in $H$. In the latter case, there is a path in $G-e$ between the two ends of $e$, thus $c((G-e)-S)=c(G-S)$. So we may suppose that $S$ is contained in $T_4$ and $\frac{|S|}{c(G-e)}<\tau(G)$. However, it is easy to see that if we contract $H$ into one vertex ($t'_1$) to obtain $G'$, then  $\frac{|S|}{c(G'-e)}<\tau(G')$ holds as well. This is a contradiction since $G'-e$ is a necklace graph, so it cannot contain such a cutset $S$.
  
   Now suppose that $r_4$ has a sibling joined to $s_1$; possibly, there are further siblings joined to that. Let $H$ denote the subgraph spanned by all these siblings. In this case, $s_1$ is a cut-vertex. Let $S=\{s_1,s_2,\ldots,s_{\ell+1}\}$. Now $G-S$ consists of all middle vertices of the bracelet $B_{\ell}$, the isolated vertices $w$ and $t'_1$, plus the nonempty component $H-s_1$. This implies that $\frac{|S|}{c(G-S)}=\frac{\ell+1}{2\ell+3}<\frac{1}{2}$, a contradiction.
   
      We obtain the same contradiction if $r_4$ has a sibling joined to both $s_1$ and $t'_1$ (i.e., to both ends of the substructure).
      
  \textit{Subsubcase (i:$\beta$): The parent $r_5$ of $r_4$ is not the root.} Thus, the parent of $r_5$ is a parallel node, which means that there is a path connecting $s_1$ and $t'_1$ that does not contain any other vertices of  $T_4$ and $T'_4$.
  
  Let $e=s_1w$, suppose that $\tau(G-e)<\tau(G)$ and let $S$ be a minimum-size tough set. By Claim~\ref{claim} we know that $s_1, w\notin S$. 
  
  First, consider the case when $s_{\ell+1}\notin S$. We can apply Lemma~\ref{l3} with $u=s_1$ and $v=s_{\ell+1}$, which implies that $S$ is either fully contained in $T_4$ or is fully outside of $T_4$. In either case, $G-e$ contains a path between $s_1$ and $w$, a contradiction.
  
  Thus, we can assume that $s_{\ell+1}\in S$. Suppose also that $s_i\in S$ and $s_j\notin S$ for some $1<i<i+1<j<\ell+1$, this contradicts Lemma~\ref{l3} with $u=s_1$ and $v=s_{j}$. Thus for some $1\le i\le\ell+1$ we have $s_1,\ldots,s_{i-1}\notin S$ and $s_i,\ldots,s_{\ell+1}\in S$. However, this $S$ cannot be minimal if $i<\ell+1$, because if we remove $s_i,\ldots,s_{\ell}$ from $S$ then the size of the cutset decreases by $\ell-i+1$ and the number of components decreases by $2(\ell-i+1)$, so the ratio (that was less than $\frac{1}{2}$) decreases. (It is easy to see that the smaller set is still a cutset.) 
  Therefore $s_1,\ldots,s_{\ell}\notin S$. This implies that in $G-S$ the vertices $w, s_1,\ldots,s_{\ell}$ and the two middle vertices between $s_{\ell}$ and $s_{\ell+1}$ (denote them by $v_1$ and $v_2$) are all in one component. 
  
  By Lemma~\ref{k}, $S$ is a tough set in $G$ as well, so $\frac{|S|}{c(G-S)}=\frac{1}{2}$. Let $S'=S-s_{\ell+1}$. Clearly $|S'|=|S|-1$. The vertex $s_{\ell+1}$ has 4 neighbors: $w, v_1, v_2, t'_1$. Since the first 3 are in one component of $G-S$, we have $c(G-S')\ge c(G-S)-1$. Using the Mediant Inequality it can be verified that $\frac{|S'|}{c(G-S')}<\frac{1}{2}$, a contradiction.

  \medskip
  \textit{Subcase (ii): $T'_4$ is $R_2$.}
   If there is no other vertex or edge in the graph, then let $e=ws_{\ell+1}$. Thus $G-e$ is nearly a necklace graph, the last edge in the series connection is missing here. However,   its toughness is still $\frac{1}{2}$, as can be shown using the argument from the proof of Lemma~\ref{necklace}. This again gives a contradiction.
   
  If $G$ also contains the leap edge $f=s_1t'_1$ then $\tau(G-f)\ge \tau(G-f-e)\ge \frac{1}{2}$, a contradiction again.

  \textit{Subsubcase (ii:$\alpha$): The parent $r_5$ of $r_4$ is the root.} We obtain a contradiction in essentially the same way as Subsubcase (i:$\alpha$), using again that the toughness of the graph obtained from a necklace graph by deleting the last edge is still $\frac{1}{2}$.
  
  \textit{Subsubcase (ii:$\beta$): The parent $r_5$ of $r_4$ is not the root.} We obtain a contradiction in essentially the same way as Subsubcase (i:$\beta$).

  \medskip
  \textit{Subcase (iii): $T'_4$ is $R_3$.}  If there is no other vertex or edge in the graph, then let $e=ws_{\ell+1}$. Thus, $G-e$ is still very similar to a necklace graph.  
  With a slight adjustment to the argument of Lemma~\ref{necklace} one can prove that 
   its toughness is still $\frac{1}{2}$, a contradiction.
  
  If $G$ also contains the leap edge $f=s_1t'_1$ then $\tau(G-f)\ge \tau(G-f-e)= \frac{1}{2}$, a contradiction again.
  
  Now suppose that $G$ contains a vertex $z$ outside of $T_4$ and $T'_4$. 
    Let $S=\{s_1,\ldots, s_{\ell+1},t'_1\}$. Thus, in $G-S$ all the middle vertices of the bracelets become isolated vertices, as well as the 3 middle vertices of $T'_4=R_3$, and $w$. There is at least one additional component that contains $z$. Therefore,  $\tau(G)\leq \frac{\ell+2}{2\ell+5}<\frac{1}{2}$, a contradiction.

  \medskip
  \textit{Subcase (iv): $T'_4$ is a parallel join of $B_{k}$ and $P_3$.}
  Let $\{s_1,\ldots, s_{\ell+1}\}$ and $\{s'_1,\ldots, s'_{k+1}\}$ denote the joining vertices of $B_{\ell}$ and  $B_{k}$ respectively so that $s'_1=s_{\ell+1}$, and  let $w$ and $w'$ denote the middle vertices of the two $P_3$ substructures.
  
  If there is no other vertex or edge in the graph, then let $e=ws_{\ell+1}$ and $f=s_{\ell+1}w'$. Now $G-e-f$ is a necklace graph, which implies that $\tau(G-e)\ge \tau(G-e-f)=\frac{1}{2}$, a contradiction again. 
  
  If $G$ also contains the leap edge $g=s_1s'_{k+1}$ then $\tau(G-g)\ge \tau(G-f-e-g)= \frac{1}{2}$, a similar contradiction.

  Thus, we can assume that $G$ contains a vertex $z$ outside of $T_4$ and $T'_4$. 
   Let $S=\{s_1,\ldots, s_{\ell+1}\}\cup \{s'_1,\ldots, s'_{k+1}\}$. Now in $G-S$ all the middle vertices of the bracelets become isolated vertices, as well as $w$ and $w'$.  There is at least one more component that contains $z$. So,  $\tau(G)\leq \frac{\ell+k+1}{2\ell+2k+3}<\frac{1}{2}$, a contradiction.
  
  \subsubsection*{Case 4(b):  $r_4$ has only two children $B_{\ell}$ and $P_3$, and $r_4$ is the root.}
  
  This graph is nearly identical to the one in Case 1(b), except that in this case, a \( P_3 \) replaces an edge. Nevertheless, we can apply essentially the same argument here to establish that \( \tau(G) > \frac{1}{2} \).
  \medskip
  
  This completes the proof that no subtree can have height 4 or greater.
  
  To summarize, we ruled out the possibility of the sp-tree having height 0, 2, or 4 above, and 
  determined the only possible graphs when the sp-tree has height 1 or 3. This concludes the proof of Lemma~\ref{heights}.
   \end{proof}
   
  Since $P_3$ is a necklace graph, Lemma~\ref{heights} implies that a reduced sp-graph is minimally $\frac{1}{2}$-tough if and only if it is a necklace graph. 
  It is easy to describe the graphs whose reduced graph is a
  necklace graph: they can be obtained from necklace graphs
  by repeatedly replacing a $P_3$ subgraph whose middle
  vertex has degree 2 with a path of length at least 3.

 We define a \emph{pearl} as the parallel join of two paths, each of length at least 2. Now the above Lemmas and Theorems imply the following characterization of minimally $\frac{1}{2}$-tough sp-graphs, which is the main result of our paper.
  
\begin{figure}[ht]
	\centering
	\scalebox{1}{
%
%
%

	\begin{tikzpicture}[scale=1,
	every edge/.style = {draw,very thick},
	vertex/.style args = {#1 #2}{circle, 
		draw, 
		thick,
		minimum size=2pt,inner sep=-2pt,
		fill=black,
		label=#1:#2},
			leaf/.style={rectangle, draw=black, rounded corners=1mm, 
			text centered, anchor=north, text=black,text width=1em,inner sep=1.5pt},
		parallel/.style={circle, draw=black, fill=gray!30,
			text centered, anchor=north, text=black,text width=1em,inner sep=1pt},
		non/.style={circle, 
				text centered, anchor=north, text=black,text width=1em,inner sep=1pt},
		series/.style={circle, draw=black, fill=black!70, 
			text centered, anchor=north, text=white, text width=1em,inner sep=1pt},
		level distance=0.25cm, growth parent anchor=south,level/.style={sibling distance=10mm}]
		
	
	\begin{scope}[shift=({0,0})]
		\foreach \x in {1,4,6,10}{
		\begin{scope}[shift=({\x,0})]
			\path   
			node (A)  [vertex=above $$]    at ( 0, 0) {} 
			node (B)  [vertex=above $$]   at ( 1, .5) {} 
			node (C) [vertex=above $$]     at ( 2, 0) {} 
			node (D) [vertex=below $$]     at ( 1, -.5) {}

			(A) edge (B)
			(B) edge (C)
			(A) edge (D)
			(D) edge (C)
			
			;
		\end{scope}
	}
		\foreach \x in {0,3,8,9,12}{
		\begin{scope}[shift=({\x,0})]
			\path   
			node (A)  [vertex=above $$]    at ( 0, 0) {} 
			node (B)  [vertex=above $$]   at ( 1, 0) {}

			(A) edge (B)

			;
		\end{scope}
	}
		\foreach \x in {1.5,4.5,5.5}{
		\begin{scope}[shift=({\x,.25})]
			\path   
			node (A)  [vertex=above $$]    at ( 0, 0) {} 
			;
		\end{scope}
	}
			\foreach \x in {1.5,2.5,6.5}{
		\begin{scope}[shift=({\x,-.25})]
			\path   
			node (A)  [vertex=above $$]    at ( 0, 0) {} 
			;
		\end{scope}
	}
	
				\foreach \x in {.5,8.5,9.5}{
		\begin{scope}[shift=({\x,0})]
			\path   
			node (A)  [vertex=above $$]    at ( 0, 0) {} 
			;
		\end{scope}
	}
		
	\end{scope}

\end{tikzpicture}

}
	\caption{A minimally $\frac{1}{2}$-tough graph}
	\label{fig:necklace2}
\end{figure}
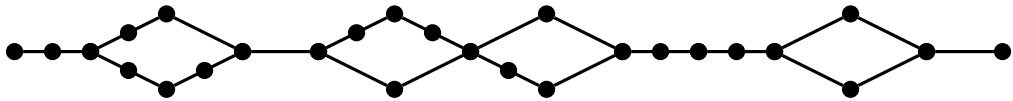

 \begin{thm}\label{mainthm}
 	An sp-graph $G$ is minimally $\frac{1}{2}$-tough if and only if it is a series connection of an arbitrary number of pearls and at least two edges, with the first and last components being edges. (See Fig.~\ref{fig:necklace2}.)
 \end{thm}
 
 Every necklace graph contains a vertex of degree 1.
 Moreover, if $\tau(G)>\frac12$, then $G$ cannot have a
 vertex of degree 1, since its only neighbor would be a
 cut-vertex. On the other hand, every sp-graph contains
 a vertex of degree at most 2. The Generalized Kriesell
 conjecture therefore follows for minimally $t$-tough
 sp-graphs with $t\geq\frac12$.
 
 \begin{thm}
 	The Generalized Kriesell conjecture holds for minimally $t$-tough series--parallel graphs if $t\ge \frac{1}{2}$.
 \end{thm}
  
  \section{Open questions}\label{open}

  The most immediate open problem is whether the Generalized Kriesell conjecture holds for minimally \(t\)-tough graphs when \(t<1/2\). A characterization of these graphs is another important open problem.
   Our methods may be applicable when \( t = \frac{1}{k} \) for some integer \( k \geq 3 \), as Lemma~\ref{k} can be used in these cases. However, we do not yet see how to characterize these graphs.
   
   Another open problem is the characterization of minimally tough \emph{generalized sp-graphs}. In these graphs, in addition to series and parallel join, one can also use the \emph{source merge} operation, i.e., identify the sources of two smaller generalized sp-graphs.
  
  
  \bibliographystyle{amsplain}

\begin{thebibliography}{xbib}
	\bibitem{t3}
	H. Broersma, E. Engbers, and H. Trommel,
	\newblock Various results on the toughness of graphs,
	\newblock \textit{Networks},
	\newblock \textbf{33} (1999) 233--238.

	\bibitem{toughness_intro}
	V. Chv\'{a}tal,
	\newblock Tough graphs and Hamiltonian circuits,
	\newblock \textit{Discrete Mathematics},
	\newblock \textbf{5}(3) (1973) 215--228.

\bibitem{mediant}
	N. Chuquet,
	\newblock \textit{Le Triparty en la Science des Nombres} (1484).

	\bibitem{t5}
	F. Gavril,
	\newblock The intersection graphs of subtrees in trees are exactly the chordal graphs,
	\newblock \textit{Journal of Combinatorial Theory Series B},
	\newblock \textbf{16}(1) (1974) 47--56,
	\newblock \url{doi:10.1016/0095-8956(74)90094-x}

	\bibitem{kriesell}
	T. Kaiser,
	\newblock Problems from the workshop on dominating cycles,
	\newblock \url{http://iti.zcu.cz/history/2003/Hajek/problems/hajek-problems.ps} (2003)

	\bibitem{specgraph}
	G.Y. Katona and K. Varga,
	\newblock Minimal toughness in special graph classes,
	\newblock \textit{Discrete Mathematics \& Theoretical Computer Science},
	\newblock \textbf{25}:3 special issue ICGT'22,
	\newblock \url{https://doi.org/10.46298/dmtcs.10180}

	\bibitem{t2}
	G.Y. Katona, D. Soltész, and K. Varga,
	\newblock Properties of minimally t-tough graphs,
	\newblock \textit{Discrete Mathematics},
	\newblock \textbf{341} (2018) 221--231.

	\bibitem{kawano}
	S. Kawano and S. Nakano,
	\newblock Generating All series--parallel Graphs,
	\newblock \textit{IEICE Transactions on Fundamentals of Electronics, Communications and Computer Sciences},
	\newblock \textbf{E88-A}:5 (2005) 1129--1135.

	\bibitem{t4}
	H. Ma, X. Hu, and W. Yang,
	\newblock On the minimum degree of minimally $1$-tough, triangle-free graphs and minimally $\frac{3}{2}$-tough, claw-free graphs,
	\newblock \textit{Discrete Mathematics},
	\newblock \textbf{346} (2023) 113352.

	\bibitem{mader}
	W. Mader,
	\newblock Eine Eigenschaft der Atome endlicher Graphen,
	\newblock \textit{Archiv der Mathematik},
	\newblock \textbf{22} (1971) 333--336.

	\bibitem{mediantproof}
	R. B. Nelsen,
	\newblock Proof without Words: Regle des Nombres Moyens,
	\newblock \textit{Mathematics Magazine},
	\newblock \textbf{67}(1) (1994) 34.

	

	\bibitem{zheng-sun}
	W. Zheng and L. Sun,
	\newblock Disproof of a conjecture on minimally $t$-tough graphs,
	\newblock \textit{Discrete Mathematics},
	\newblock \textbf{347} (2024) 113982.
\end{thebibliography}

 	\end{document}